\documentclass[11pt,onecolumn]{IEEEtran}

\usepackage{amsmath,amssymb} 
\usepackage{dsfont}
\usepackage{graphicx}
\usepackage{color}
\usepackage{bm}
\usepackage{algorithmic}
\usepackage{algorithm}
\usepackage{todonotes}
\usepackage{caption}
\usepackage{subcaption}
\usepackage{tablefootnote}

\begin{document}

\title{An optimization-based approach to automated design}


\author{Ion Matei, Maksym Zhenirovskyy, John Maxwell and Johan de Kleer
}

\maketitle

\begin{abstract}
We propose a model-based, automated, bottom-up approach for design, which is applicable to various physical domains, but in this work we focus on the electrical domain. This bottom-up approach is based on a meta-topology in which each link is described by a universal component that can be instantiated as basic components (e.g., resistors, capacitors) or combinations of basic components via discrete switches. To address the combinatorial explosion often present in mixed-integer optimization problems, we present two algorithms. In the first algorithm, we convert the discrete switches into continuous switches that are physically realizable and formulate a parameter optimization problem that learns the component and switch parameters while inducing design sparsity through an $L_1$ regularization term. The second algorithm uses a genetic-like approach with selection and mutation steps guided by ranking of requirements costs, combined with continuous optimization for generating optimal parameters. We improve the time complexity of the optimization problem in both algorithms by reconstructing the model when components become redundant and by simplifying topologies through collapsing components and removing disconnected ones. To demonstrate the efficacy of these algorithms, we apply them to the design of various electrical circuits.
\end{abstract}

\section{Introduction}
\begin{figure}[!htp]
\centering
\includegraphics[width=15cm]{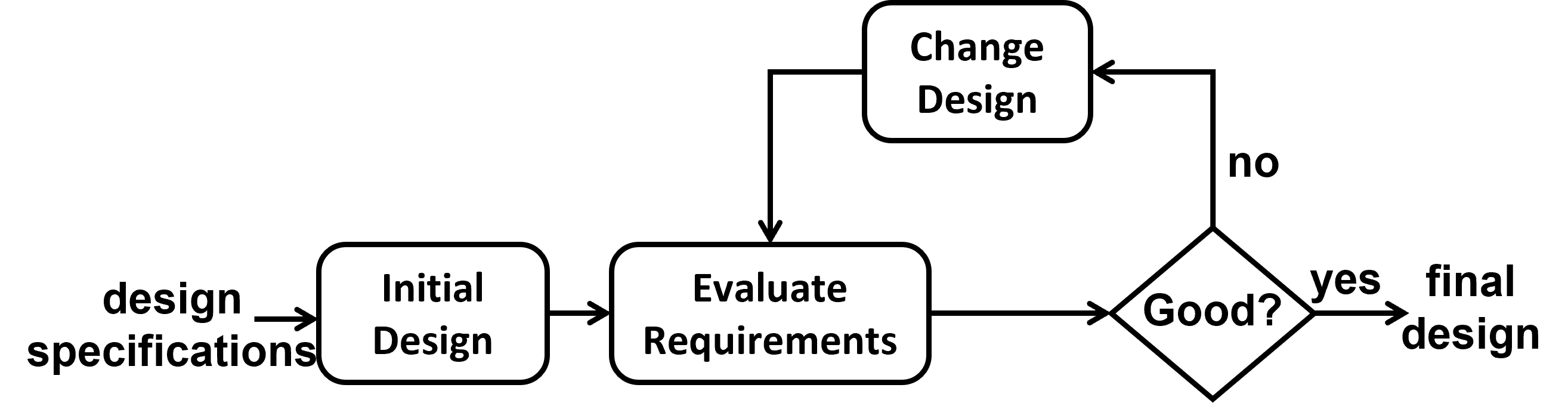}
\caption{The design optimization process: an initial design is continuously refined until requirements are satisfied.}
\label{fig:design_process}
\end{figure}

A simplified view of the design process is depicted in Figure \ref{fig:design_process}. A designer specifies a set of design requirements and specifications (e.g., components, interfaces) and begins with an initial design, which is continuously refined until it satisfies the requirements.
There are two main approaches for implementing the design process: top-down and bottom-up. In the top-down design approach, the system is first divided into subsystems without providing details about their internals. This is followed by several refinements where more information about the subsystem internals is provided, and additional subsystem levels may be added. This process continues until the description of the system is reduced to basic components. In the bottom-up approach, the system starts with basic components that are combined into subsystems with increasing complexity, which are then further combined until a system with the desired requirements emerges.

In this paper, we describe a general approach for designing physical systems using a bottom-up approach that implements the ``change design'' process in Figure \ref{fig:design_process}. This type of problem can be formulated as a mixed integer program that includes a combinatorial part to  select the component types and a continuous optimization part that selects parameters of components to meet requirements. A brute force approach to solving such an optimization problem suffers from combinatorial explosion, and heuristics based on branch-and-bound methods do not scale with the number of discrete optimization variables \cite{Clausen2003BranchAB,MORRISON201679}. To limit the effects of combinatorial explosion, we introduce two algorithms in this paper.

To facilitate the description of the algorithms and results, we focus on a design problem in the electrical domain, where the goal is to build systems with various filtering properties. However, the approach can be easily extended to other physical domains. We use the Modelica language to describe the basic components and the generated design solutions, which allows subject matter experts to interpret and evaluate the generated designs.

\begin{figure}[!htp]
\centering
\includegraphics[width=15cm]{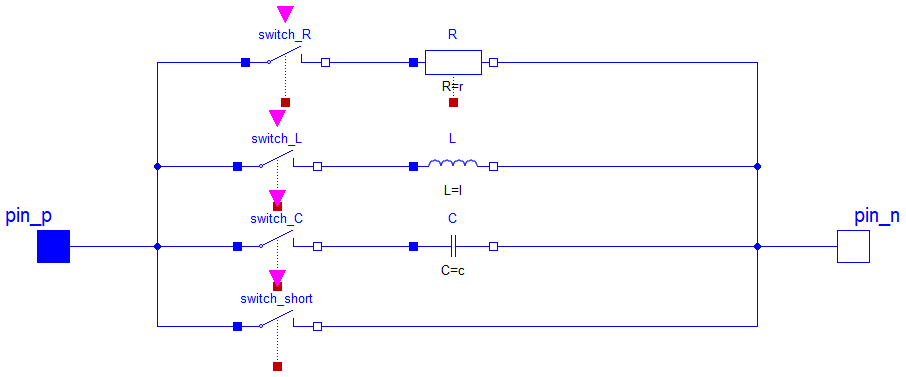}
\caption{Modelica model for the universal component. Each basic component has associated a switch to activate or deactivate it.}
\label{fig:uniniversal component}
\end{figure}

Both algorithms use a universal component that embeds the behavior of basic components in the electrical domain: resistor, inductor, capacitor, short connection, and open connection. We later show how to naturally include more complex components such as operational amplifiers (OpAmps). The schematic of the universal component is shown in Figure \ref{fig:uniniversal component}.
Each component is activated or deactivated by a switch that controls the flow of current through. The universal component has 9 modes: open, short, and all possible combinations of resistor, inductor, and capacitor components in parallel. The design problem is to find the correct switch assignments and component parameter values to meet the requirements, which can be specified in terms of the time evolution of certain quantities of interest or the characteristics of a transfer function in the case of filter design. We start with a topology, which is a diagram describing how the universal components are connected and includes points for setting boundary conditions and taking measurements. The design problem is then formulated as an optimization problem that minimizes a loss function $\mathcal{C}(\hat{\boldsymbol{y}}_{0:T}(\boldsymbol{p}, \boldsymbol{s}),\boldsymbol{y}_{0:T})$, where $\boldsymbol{p}$ and $\boldsymbol{s}$ are the parameters and switches of the basic components, respectively, $\boldsymbol{y}_{0:T}$ is a target vector of measurements over time interval $[0,T]$, $\hat{\boldsymbol{y}}_{0:T}(\boldsymbol{p}, \boldsymbol{s})$  is the model's predicted measurements, and $\mathcal{C}$ is a metric that measures the error between the model predictions and the target measurements (e.g., mean square error). The optimization problem also takes into account dynamic constraints, bounds on component parameters (e.g., resistances must be non-negative), and integer constraints (e.g., the entries of the switch vector must be 0 or 1). The main contributions of this paper are as follows::
\begin{itemize}
  \item \textit{Continuous relaxation with lossless realization:} We developed an optimization algorithm that relaxes the integer constraints on the switches by treating them as continuous variables in the range [0, 1]. The parameters of the components and their associated switches are optimized using gradient-free search algorithms and simulations based on Functional Mockup Units (FMUs) \cite{Blochwitz11thefunctional}. To encourage sparsity in the design solution, we also add an $L_1$ regularization term to the loss function. The non-zero switches are not approximated by 0 or 1, but are realized as electric resistors, ensuring no loss in optimality but a possible loss in sparsity. Since we cannot guarantee finding the global optimum, we also use parallel optimization runs with random initial conditions to generate a diverse set of design solutions.
  \item \textit{Pseudo-random search to control the combinatorial explosion:} We developed an iterative procedure that learns design solutions by sampling from the switch distributions and optimizing the parameters of the selected components. The optimization results are ranked, and the best intermediate designs are used in the next iteration, where we again sample from the surviving switch distributions. To speed up the learning process, we run the optimization algorithms in parallel. The algorithm stops when further simplification of the design through component elimination does not result in an improvement.
  \item \textit{Scalability improvement via model simplifications:} During optimization, when certain components are no longer needed (i.e., their switches are set to zero), we eliminate them and reconstruct the design model. This reduces the complexity of the model, as measured by the number of equations, and leads to faster simulation times. In addition, we developed a graph theory-based algorithm that further simplifies the designs generated by the optimization procedure. The algorithm removes unnecessary components, combines compatible components in series and parallel connections into equivalent single components, and annotates the resulting design models for visual representation and simulation in tools that support the Modelica language.
\end{itemize}

\emph{Paper structure:} In Section \ref{sec:design_optimization}, we present two algorithms for automated design: one that uses continuous relaxation, and another that combines random search with continuous optimization. In Section \ref{sec:model_simplification}, we discuss how we improve the efficiency of our design algorithms by reducing the complexity of the design models that are simulated during the design space exploration. Finally, in Section \ref{sec:results}, we present the designs generated by the two algorithms for various circuit design problems, and show how the design search and simplification process evolves until reaching a satisfactory design solution.

\section{Design optimization}
\label{sec:design_optimization}
In this section, we present two approaches for addressing the combinatorial explosion often present in design problems.
\subsection{Continuous relaxation with lossless realizations}

When using branch-and-bound heuristics to solve mixed integer programs, we may encounter situations where the cost of the relaxed problem is better than the cost obtained by converting the optimization variables into integer values. In this section, we present a method to avoid such a case. The key idea is to interpret the switches in a way that allows for their physical implementation, even when they do not have integer values. In the universal component definition, each basic component has a corresponding switch that opens or closes a connection. When the switch is open, no current flows through the component, leading to its exclusion from the design model. In the Modelica electrical library, the switch ({\tt Modelica.Electrical.Analog.Ideal.IdealOpeningSwitch}) is modeled such that when the switch is open, there is a high resistance that blocks the flow of current. When the switch is closed, there is a very low resistance that allows the current to flow freely. This switch is controlled by a boolean input, the value of which determines the switch mode. We draw inspiration from the definition of the Modelica switch to create a continuous switch that is controlled by a parameter that takes values in the range [0, 1]. The switch is defined by the equations
\begin{eqnarray}
\label{eq:11021606}
  v &=& a((\varepsilon-1)s+1), \\
  \label{eq:11021608}
  i &=& a((1-\varepsilon)s+\varepsilon),
\end{eqnarray}
where $v$ is the switch voltage, $i$ is the current through the switch, $a$ is an auxiliary variable, $s\in [0,1]$ is the switch control, and $\varepsilon$ is a small hyper-parameter that determines the residual resistance when the switch is closed. The switch equation can be simplified to
$$v = \frac{(\varepsilon-1)s+1}{(1-\varepsilon)s+\varepsilon}i,$$
showing that for $s=0$ we have $v = i/{\varepsilon}$ and for $s=1$ we have $v = \varepsilon i$, the expected behavior of a switch. We do not use this simplified representation of the switch for numerical stability reasons. The introduction of the auxiliary variable $a$ prevents the presence of equations with terms that involve divisions by very small numbers.  However, the disadvantage is that the resulting system of equations for the design model becomes a differential algebraic equation (DAE) rather than an ordinary differential equation (ODE). This limitation restricts the type of optimization approach that can be used, as we cannot directly utilize platforms that support automatic differentiation (AD) (see Section \ref{sec:discussion}). In addition to the requirements loss function $\mathcal{C}$,  we introduce a sparsity-promoting $L_1$ regularization term, resulting in the total optimization loss:
$$\mathcal{L}(\boldsymbol{p}, \boldsymbol{s}) = \mathcal{C}(\hat{\boldsymbol{y}}_{0:T}(\boldsymbol{p}, \boldsymbol{s}),\boldsymbol{y}_{0:T}) + \lambda \|\boldsymbol{s}\|_1,$$
where $0\leq {s}_i\leq 1$, with $\boldsymbol{s} = (s_i)$, and $\lambda$ is a positive weight that controls the sparsity strength. If in the optimization solution not all entries of $\boldsymbol{s}$ are zero or one, we map them into electric resistors with equivalent resistances, $\frac{(\varepsilon-1)s_i+1}{(1-\varepsilon)s_i+\varepsilon}$. Thus, we can physically realize them, without affecting the optimization cost function, i.e., the design requirements.

The pseudocode for this algorithm is shown in Algorithm \ref{alg:1}. We use a gradual approach to achieve sparsity. We start with a small $\lambda$ value to make sure that we generate an initial design that satisfies the requirements. Then we gradually increase $\lambda$ until the requirements cost function is no longer improved. Ideally, for each $\lambda$, we would like to obtain the optimal solution.  The strategy for updating $\lambda$ is reminiscent to a primal-dual approach \cite{Bertsekas_99}, where we minimize $\mathcal{C}$ under an $L_1$ sparsity constraint. In the discussion section, we make a more clear parallel between Algorithm \ref{alg:1} and a constrained optimization approach to solve multi-objective optimization problems.

In our approach, we incrementally increase the value of $\lambda$ until it begins to negatively impact the requirements cost function. At this point, we halt the process and perform a final optimization without the $L_1$ regularization term. The result of this final optimization will be our chosen design solution. Box constraints are commonly used in our problem setup, but we  use variable transformations to eliminate them and use an unconstrained optimization algorithm to minimize $\mathcal{L}$. For example, we can eliminate the constraint $a\leq x\leq b$ by using the transformation $x = a + (\sin(\tilde{x}) + 1)(b-a)/2$, where $\tilde{x}$ is the new optimization variable. It is not guaranteed that the optimization will converge to the global minimum, as the cost function's nonlinear dependence on the optimization parameters means we cannot accurately predict the structure of the problem. Ideally, we would find at least a local minimum for each $\lambda$ value, but it is possible that the optimization algorithm may take too many iterations to converge. As a result, we set a limit on the number of iterations allowed between $\lambda$ updates for practical reasons.

All optimization algorithms will require the evaluation of the design model. We use a black-box approach to optimization, where the cost evaluation is done by querying a computational model of the design: an FMU \cite{Blochwitz11thefunctional}. In the cosimulation version of the FMU, such a representation contains the algorithm used for simulating the model (e.g., CVODE solver \cite{hindmarsh2005sundials}), in addition to the design description. FMUs can be integrated in several languages (e.g., Python, C, Java) and computational platforms (e.g., Matlab/Simulink, OpenModelica, Dymola). We used a gradient free (i.e., a direct method) optimization algorithm that relies only on the objective function, namely Powell's method \cite{10.1093/comjnl/7.2.155}. Empirically, it provides an better convergence rate than other gradient-free algorithms such as Nelder-Mead.

\begin{algorithm}
\caption{Continuous relaxation design algorithm}
\label{alg:1}
\begin{algorithmic}[1]
\REQUIRE $\delta$: solution tolerance
\REQUIRE $\lambda$: $L_1$ loss weight
\REQUIRE $\Delta$: $L_1$ loss weight increase rate
\REQUIRE FMU of the initial design model
\REQUIRE $\boldsymbol{p}$, $\boldsymbol{s}$: initial parameter and switch values
\REQUIRE $\boldsymbol{y}_{0:T}$: target measurements
\STATE {$\mathcal{C}_{prev}=\infty$}
\WHILE{True}
\STATE $\boldsymbol{p}, \boldsymbol{s}\leftarrow \arg\min_{\boldsymbol{p}, \boldsymbol{s}} \mathcal{C}(\hat{\boldsymbol{y}}_{0:T}(\boldsymbol{p}, \boldsymbol{s}),\boldsymbol{y}_{0:T})+\lambda \|\boldsymbol{s}\|_1$
\STATE {$\mathcal{C}^*=\mathcal{C}(\hat{\boldsymbol{y}}_{0:T}(\boldsymbol{p}, \boldsymbol{s}),\boldsymbol{y}_{0:T})$}
\IF {$\mathcal{C}^* \leq \mathcal{C}_{prev}$}
\STATE $\lambda \leftarrow \Delta\lambda$
\STATE $\mathcal{C}_{prev} = \mathcal{C}^*$
\STATE eliminate components corresponding to zero switches and reconstruct the model
\ELSE
\STATE $\boldsymbol{p}, \boldsymbol{s}\leftarrow \arg\min_{\boldsymbol{p}, \boldsymbol{s}} \mathcal{C}(\hat{\boldsymbol{y}}_{0:T}(\boldsymbol{p}, \boldsymbol{s}),\boldsymbol{y}_{0:T})$
\RETURN $\boldsymbol{p}, \boldsymbol{s}$
\ENDIF
\ENDWHILE
\end{algorithmic}
\end{algorithm}

\subsection{Pseudo-random search}
An alternative to continuously relaxing the switches is to consider them as integer values from the start. To reduce the combinatorial explosion, we can limit each universal component to one of five modes: open, short, resistor, induction, or capacitor. This can be achieved by declaring an integer input for the universal component that determines the mode using Modelica {\tt if-then} statements. For example, if the input is 0, all switches on the branches of the universal component would be set to open.
Applying a branch-and-bound method to solve the mixed integer program would not be feasible because real values would not be relevant in this case. If we had to revert to continuous relaxation of the switches, the complexity would be too great due to the large number of branches that would need to be considered.

One potential solution to the combinatorial problem of selecting a topology is to use reinforcement learning (RL) \cite{Sutton1998}. For example, \cite{10.5555/3437539.3437740} used RL to solve design problems, where the state corresponds to a particular topology and current component parameter values, and the policy determines what changes to the topology are necessary. However, our own attempts to use this approach encountered insurmountable challenges: (i) the number of possible topologies was significantly larger than the one considered in \cite{10.5555/3437539.3437740}, and (ii) the time required to reach an optimal policy was too long due to the large number of iterations of the RL algorithm and the time needed to simulate the explored designs. Therefore, we decided to use a more practical approach that leveraged the parallelism of many model simulations, combined with random selections of universal component modes and component parameter optimization for a limited number of explored design models.

The random search approach is described in Algorithm \ref{alg:2},
which has three important hyperparameters: $n_s$ representing the budget for simulations,
$n_o$ the budget for optimization, and $C_{th}$ the acceptable threshold with
respect to the requirement cost. The algorithm starts by generating $n_s$ feasible
(i.e., can be simulated) topologies using as starting point a meta-design model
(line 2 of the algorithm). These topologies are generated by randomly sampling
the modes of the universal components, and their parameters (i.e., the resistance,
capacitance or inductance). All these $n_s$ topologies are simulated in parallel (line 4)
and the topologies  are ranked based on their requirement cost functions (line 5).
We choose the best $n_o$ topologies as dictated by their corresponding requirements cost
functions (line 6). By best we mean the topologies with the smallest requirements
cost functions. For each such topology we optimize the parameters of the components
by minimizing the requirements cost functions in terms of the component parameters
(line 7). Next, we check if there exists a least one optimized topology that
generates an optimal cost that is smaller than the threshold $C_{th}$.
If we found such topology, it means that we have a design that meets
the requirements. However, it does not mean that we found the most parsimonious design.
Therefore, we take the optimized requirements and generate sparser designs by
randomly setting components that are not {\tt short} or {\tt open} to these modes
(lines 15-18). While this process may appear to generate a combinatorial explosion,
in practice this is not the case. Many of the generated topologies are not
feasible (i.e., cannot be simulated). In addition, many of the topologies
generated are not unique, since, similar to the continuous relaxation case,
we execute simplification step for each topology. Such a simplification step
eliminates dangling components, often resulting in topologies that have been
previously generated. Here we use structural similarity (i.e., with respect to topology),
and do not care about the component parameter values. However, in the unlikely case
 we generate a large number of topologies, the algorithm will pick $n_s$ of them,
 based on the requirements cost of their parents (lines 21-22).
 This process continues until newly generated, sparser topologies cannot
 meet the $C_{th}$ threshold. The parameter $n_s$ is much bigger than $n_o$
 since parallel simulations of design models are much faster than optimising
 for designs. In the results section, we provide statistical results on the random
 search process leading to a design that meets requirements.
\begin{algorithm}
\caption{Cost guided random search design algorithm}
\label{alg:2}
\begin{algorithmic}[1]
\REQUIRE FMU of the initial design model
\REQUIRE $C_{th}$: cost threshold
\REQUIRE $n_s$: max number of topologies
\REQUIRE $n_{o}$: parameter to filter topologies based on cost function value
\STATE $\mathds{S}=\emptyset$
\STATE {Generate as set $\mathds{T}$ of $n_s$ feasible topologies with randomly selected modes for its universal components}
\WHILE{True}
\STATE {Simulate topologies in $\mathds{T}$ and produce set $\mathds{C}$ that includes their requirement costs }
\STATE {Sort the set $\mathds{C}$ in increasing order}
\STATE{Generate the set $\tilde{\mathds{T}}$ that includes topologies corresponding to the first $n_{o}$ costs in the ordered set $\mathds{C}$}
\STATE {For all topologies in $\tilde{\mathds{T}}$ optimize the parameters of the universal components}
\STATE {Generate the set of costs $\tilde{\mathds{C}}$ corresponding to the optimized topologies in $\tilde{\mathds{T}}$}
\IF {$\min\tilde{\mathds{C}}>C_{th}$}
\STATE {Return $\mathds{S}$}
\ELSE
\STATE {$\mathds{S} = \{\tilde{T}\in \tilde{\mathds{T}}| \textmd{the cost of }\tilde{T} \leq C_{th}\}$}
\STATE $\mathds{T}=\emptyset$
\FOR {all topologies $T$ in ${\mathds{S}}$}
\STATE {Generate new topologies by successively selecting a universal component that is not in a {\tt short} or {\tt open} mode and setting it to the mode {\tt short} and {\tt open}}
\STATE {Eliminate from the newly generated topologies the infeasible ones}
\STATE {Simplify the topologies eliminating dangling components and components with zero resistances and capacitances}
\STATE {If the remaining topologies are not in $\mathds{T}$ add them}
\ENDFOR
\IF {$|\mathds{T}|>n_s$}
\STATE {Order the topologies in $\mathds{T}$ based on the requirement costs of their parents}
\STATE $\mathds{T}\leftarrow$\{first $n_s$ topologies in the ordered $\mathds{T}$\}
\ENDIF
\ENDIF
\ENDWHILE
\end{algorithmic}
\end{algorithm}

\section{Model Construction and Simplification}
\label{sec:model_simplification}
\begin{figure}[!htp]
\centering
\includegraphics[width=15cm]{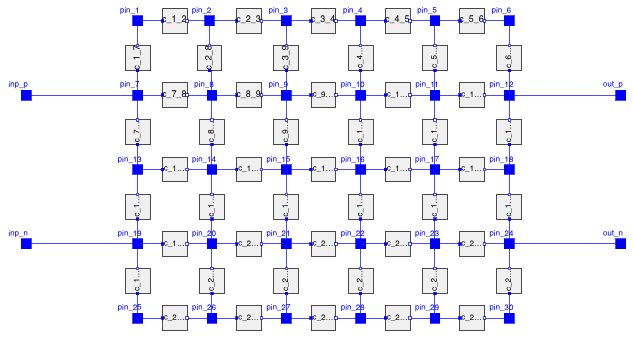}
\caption{Modelica model for the grid.  Universal components connect the grid points.}
\label{fig:5x6_grid}
\end{figure}

\begin{figure}[!htp]
\centering
\includegraphics[width=15cm]{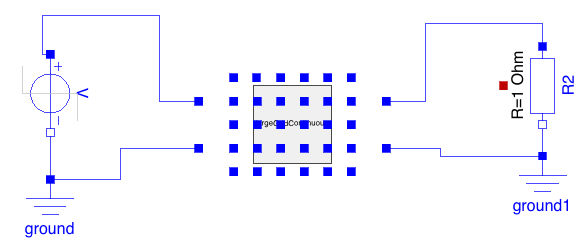}
\caption{Modelica model for the scenario that gives the boundary conditions of a grid.}
\label{fig:grid_scenario}
\end{figure}

We automatically construct a Modelica model for a domain given a universal component and a specification of the grid.  For instance, if the user wanted to use a 5x6 grid, then the program would generate a Modelica model with 30 grid points with components connecting pairs of points vertically and horizontally (see Figure \ref{fig:5x6_grid}).  This model is embedded in another model which specifies the components that set the boundary conditions, i.e., the voltage source and the resistor load (see Figure \ref{fig:grid_scenario}).

In the continuous relaxation approach to optimization, each universal component has switches that allow internal components to be enabled or disabled.  These switches can be set from the top level model.  When a component is disabled, then the Modelica compiler ignores it when constructing an FMU, thus no equations pertaining to the respective components are added. This process is implemented by conditionally declaring the basic components of the universal component. Consequently, a basic component appears in the instance of a universal component only when a corresponding flag is set to true. The flags of the basic components in all instances of the universal component are continuously updated during the optimization process.

\begin{figure}[!htp]
\centering
\includegraphics[width=15cm]{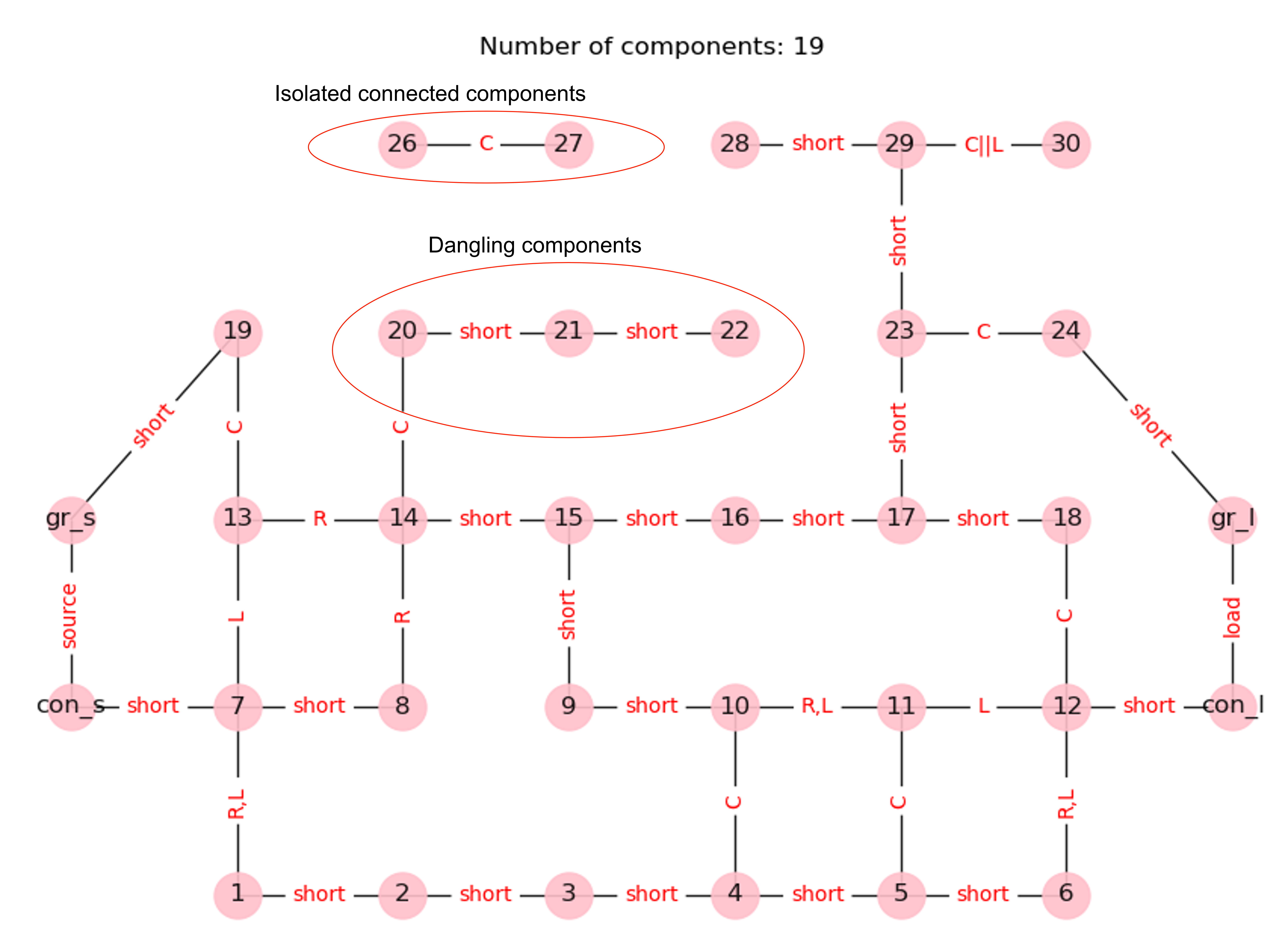}
\caption{Graph representation of a design solution: vertices are connection points and edges components.}
\label{fig:design_example}
\end{figure}

\begin{figure}[!htp]
\centering
\includegraphics[width=15cm]{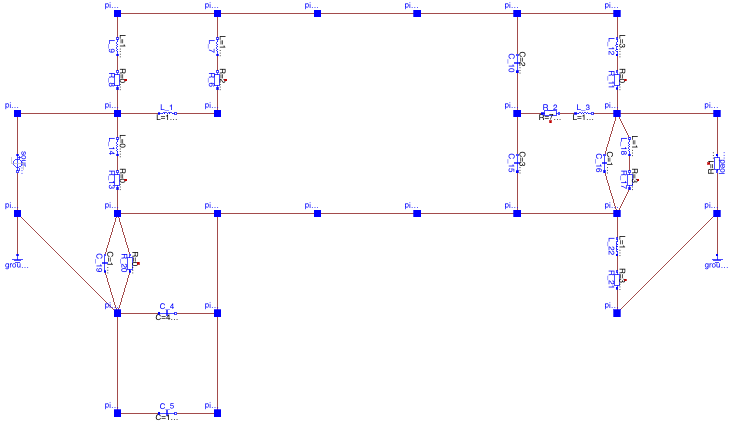}
\caption{Modelica model for a solution found by the optimizer.}
\label{fig:grid_solution}
\end{figure}

\begin{figure}[!htp]
\centering
\includegraphics[width=15cm]{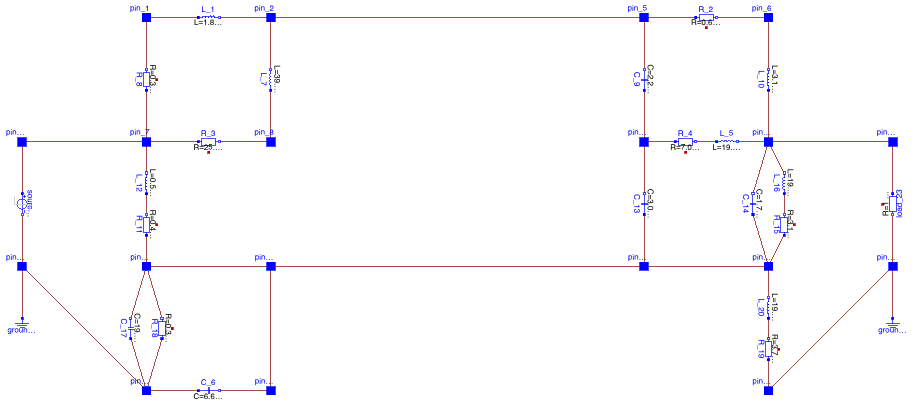}
\caption{Simplified modelica model for Figure \ref{fig:grid_solution}.}
\label{fig:simplified_grid_solution}
\end{figure}

After the optimizer has found a solution (i.e., has determined which components should be enabled and what their parameter values should be), we produce another Modelica model that flattens the universal components and just shows the internal components.  At this point we perform two simplification operations: eliminate isolated components and dangling components. These operations are necessary to deal with the cases where switch, resistor or capacitor values are close to zero. Such a situation indicates the presence of open connections. Figure \ref{fig:design_example} shows a design solution example that contains isolated (capacitor between vertices 26 and 27) and dangling (components between vertices 14,20,21,22) components.
The design solution can contain isolated components since switches are not exactly zero, meaning that there may be some very small residual currents passing true a components. Thus, it may appear that we have components that are isolated but in fact only a small, negligible current passes through them. The isolated components are eliminated by first generating the largest set of connected components that include the boundary conditions (i.e., the voltage source and the resistor load), and discarding the remaining ones. The design solution may also contain components that appear to be dangling, i.e., they are connected at one end only. The reason for such a phenomenon is the same as in the isolated components case: residual currents passing through them. The dangling components are found by looking at the cycles of the design. If a component does not belong to a cycle then it must be dangling, thus it is eliminated. Residual currents can appear even when using the pseudo-random approach, where the switches are not approximated by real numbers. The reasons are similar to the continuous relaxation case, i.e., small values of resistors and capacitors. We implemented code that generates a visually interpretable layout for the components based on the and-or graph also of the components that are between two grid points (see Figure \ref{fig:grid_solution}). The layout was achieved by annotating the flattened Modelica design model with Modelica notation that generates the visual effects.
Finally,  we have code to simplify the model by merging compatible serial or parallel components.  The code goes through this process iteratively, until no merging can be achieved. The resulting model has correct equivalent parameter values (i.e., resistances in serial connections are added) and it can be simulated using Modelica (see Figure \ref{fig:simplified_grid_solution}).

\section{Discussion}
\label{sec:discussion}
\subsection{Connections to optimization algorithms for design}
Algorithm \ref{alg:1} can be seen as a heuristic to solve a multi-objective optimization problem, aimed at jointly minimizing the requirements cost $ \mathcal{C}(\hat{\boldsymbol{y}}_{0:T}(\boldsymbol{p}, \boldsymbol{s}),\boldsymbol{y}_{0:T})$ and the sparsity cost $\|\boldsymbol{s}\|_1$. Since minimizing the requirements cost is the most important criterion, we can apply  a lexicographic method to multi-objective optimization \cite{ARORA2017771}, where we first minimize:
$$\min_{\boldsymbol{p}, \boldsymbol{s}}\mathcal{C}(\hat{\boldsymbol{y}}_{0:T}(\boldsymbol{p}, \boldsymbol{s}),\boldsymbol{y}_{0:T}),$$
to generate an optimal value $\mathcal{C}^*$. This can be achieved by starting with a very small value for $\lambda$ in Algorithm \ref{alg:1}. Next, we would solve a constrained optimization problem of the form
\begin{eqnarray*}
  \min_{\boldsymbol{p}, \boldsymbol{s}} & & \|\boldsymbol{s}\|_1 \\
  \textmd{subject to: } & & \mathcal{C}(\hat{\boldsymbol{y}}_{0:T}(\boldsymbol{p}, \boldsymbol{s}),\boldsymbol{y}_{0:T}) \leq \mathcal{C}^*,
\end{eqnarray*}
that can be solved using a penalty-based approach that converts constrained optimization problems into unconstrained problems, by adding penalty terms to the objective function: $$\min_{\boldsymbol{p}, \boldsymbol{s}} \left\{ \|\boldsymbol{s}\|_1 + \rho \left[\mathcal{C}(\hat{\boldsymbol{y}}_{0:T}(\boldsymbol{p}, \boldsymbol{s}),\boldsymbol{y}_{0:T}) - \mathcal{C}^*\right]\right\}.$$
Equivalently, we have
$$\min_{\boldsymbol{p}, \boldsymbol{s}} \left\{\lambda\|\boldsymbol{s}\|_1 + \left[\mathcal{C}(\hat{\boldsymbol{y}}_{0:T}(\boldsymbol{p}, \boldsymbol{s}),\boldsymbol{y}_{0:T}) - \mathcal{C}^*\right]\right\},$$
where $\lambda = 1/\rho$, which reflects the lines 5-7 of Algorithm \ref{alg:1}, when omitting the constant $\mathcal{C}^*$. Solving for the exact $\mathcal{C}^*$ is often impractical due to the complexity of the initial design. Consequently, we use a surrogate for $\mathcal{C}^*$, given by the cost function after a pre-set number of iterations of the optimization algorithm, where the cost function is updated after each outer iteration.

Algorithm \ref{alg:2} can be interpreted as a population-based approach to design (i.e., genetic algorithm \cite{goldberg89}), where a collection of design points is optimized.  Algorithm \ref{alg:2} includes an initialization step that generates a set of individuals (line 2), followed by a selection step, where we keep a subset of the individuals to generate new ones. The selection criterion is done by looking at their requirements cost functions (lines 4-8). Unlike  a typical genetic algorithm, we do not have explicit cross-over, but we do have a mutation step, implemented by randomly converting components (i.e., resistors, capacitors, inductors), into {\tt{short}} or {\tt{open}} connections (lines 14-18). Cross-over is not critical to be implemented since empirically we noticed that many generated topologies share some commonality, to the extent that we sometimes have to reject new topologies for having already been considered.

\subsection{Including operational amplifiers}
Integrated circuits such as OpAmps are a fundamental part of modern electronic circuits and are commonly used as building blocks in analog circuits. Therefore, it is reasonable to include them in the definition of a universal component. However, as OpAmps are three-pin components, it is not immediately clear how to incorporate them into the universal component definition. One concern is that certain configurations of OpAmps may result in designs that cannot be simulated. An easy but wasteful approach to include OpAmps is to find a design built of passive components only, and replace the inductors with OpAmp-based implementations, i.e., floating inductors. Such implementations include  two General Immittance Converters \cite{doi:10.1080/00207217508920481} and include 7 resistors, 2 capacitors and 4 OpAmps. In what follows we describe two approaches to include OpAmps that require fewer passive and active components.

\begin{figure}
\centering
\begin{subfigure}[b]{.45\textwidth}
  \centering
  \includegraphics[width=\textwidth]{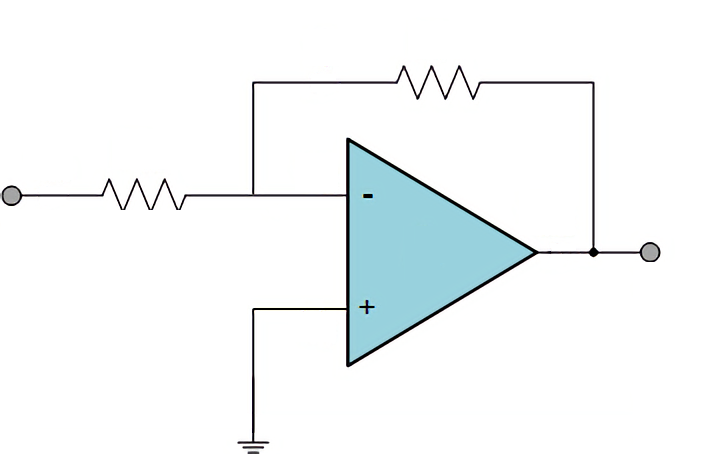}
   \caption{inverting operational amplifier.}
   \label{fig:inverting_op_amp}
\end{subfigure}
\begin{subfigure}[b]{.36\textwidth}
  \centering
  \includegraphics[width=\textwidth]{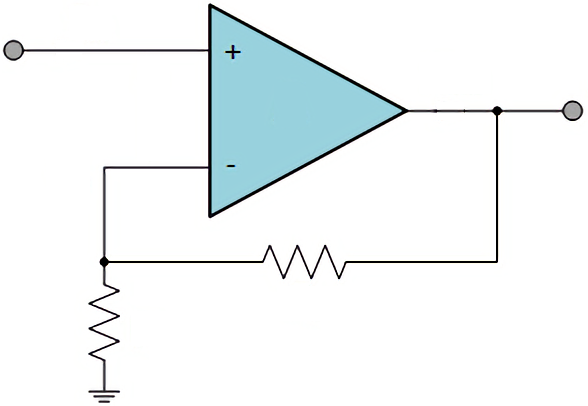}
   \caption{non-inverting operational amplifier.}
   \label{fig:non_inverting_op_amp}
\end{subfigure}
\caption{Operational amplifier configurations.}
\label{fig:op_amp_configurations}
\end{figure}

\begin{figure}
\centering
  \centering
  \includegraphics[width=15cm]{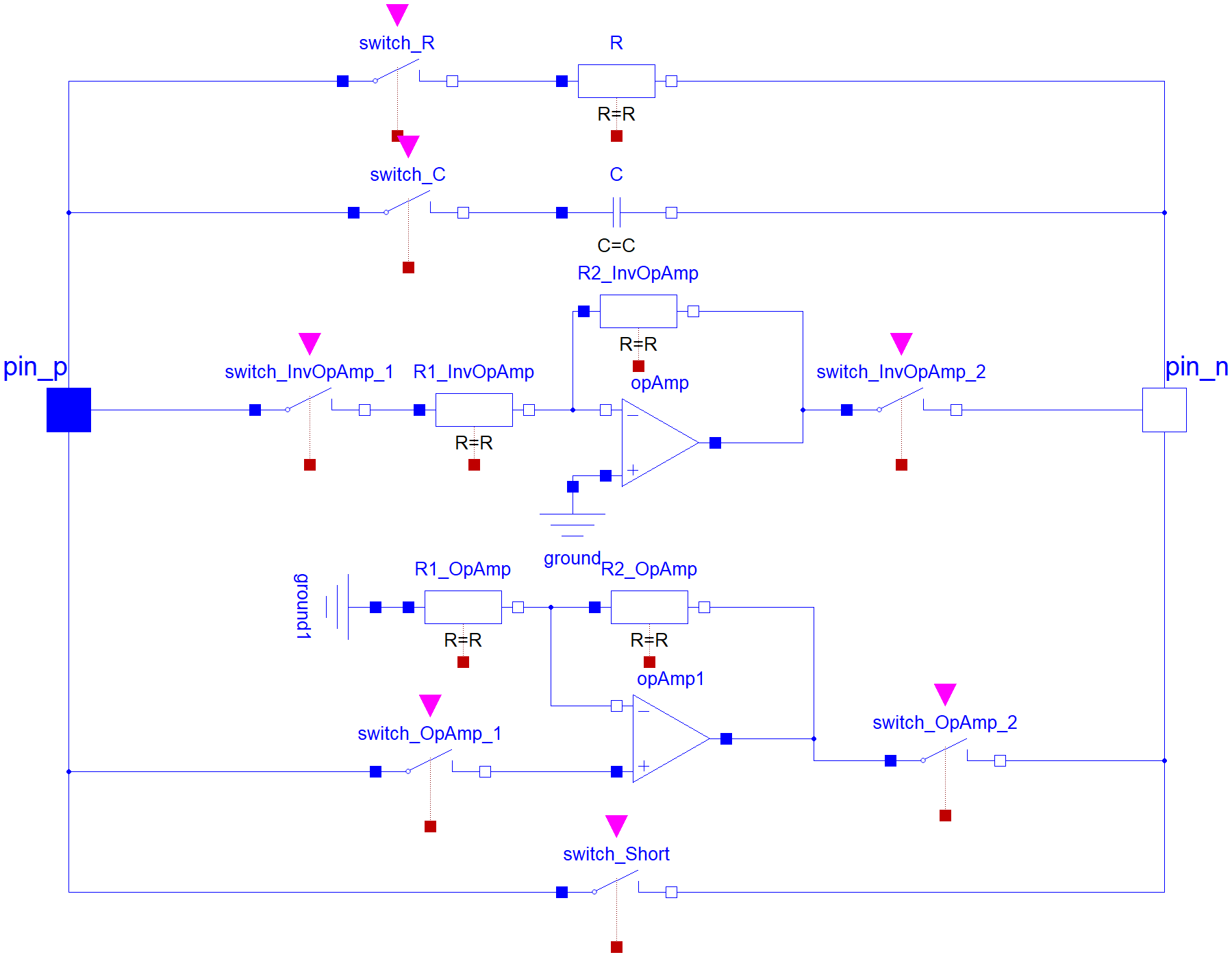}
   \caption{Universal component based on inverting and non-inverting OpAmp configurations.}
   \label{fig:universal_component_op_amps}
\end{figure}

\begin{figure}
\centering
\begin{subfigure}[b]{.45\textwidth}
  \centering
  \includegraphics[width=\textwidth]{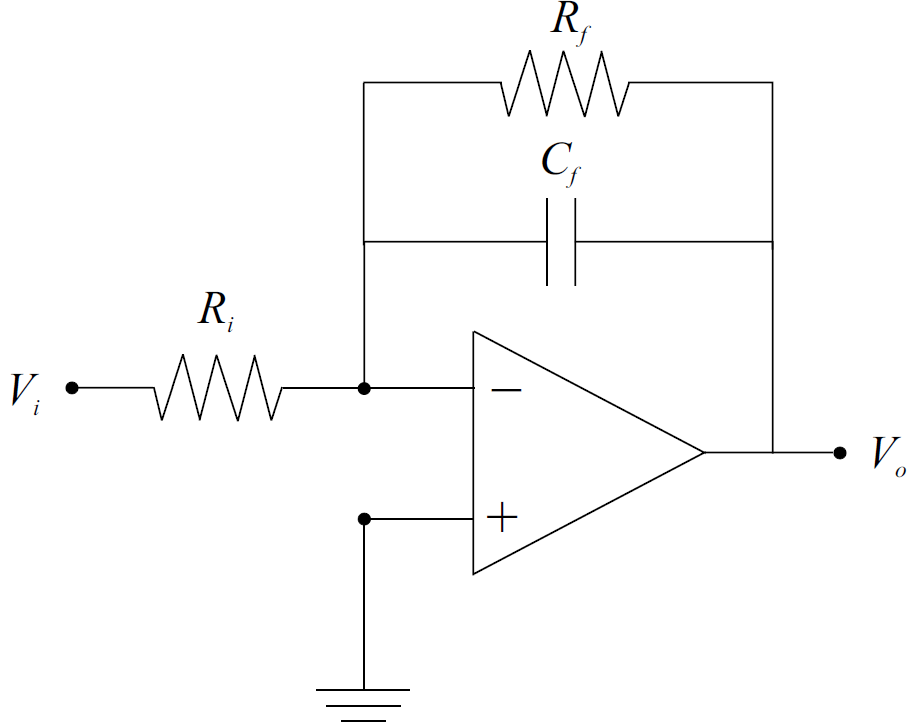}
   \caption{first order low pass filter.}
   \label{fig:folp_filter}
\end{subfigure}
\begin{subfigure}[b]{.45\textwidth}
  \centering
  \includegraphics[width=\textwidth]{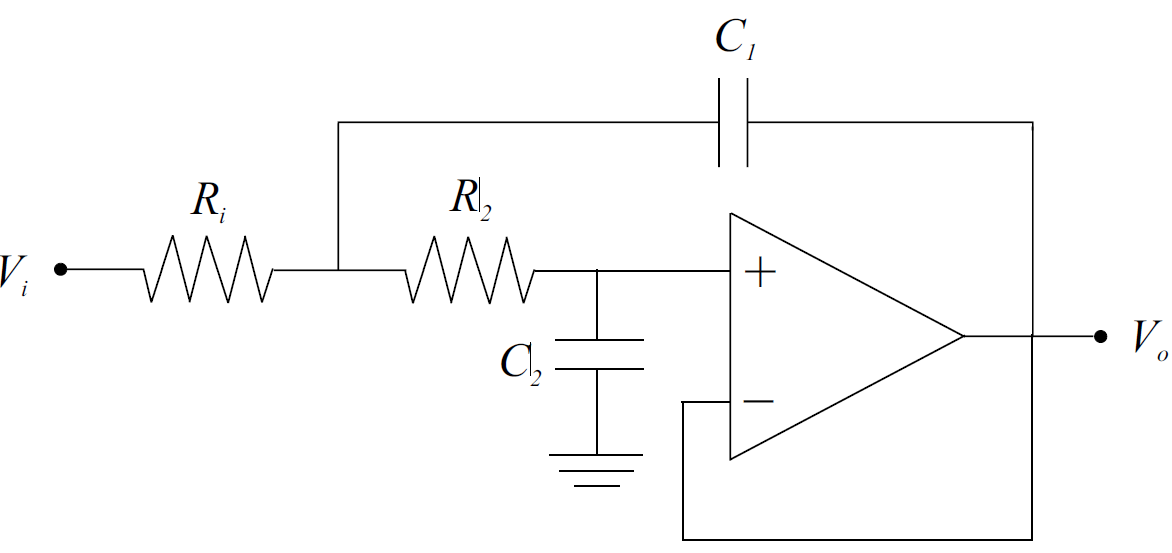}
   \caption{second order low pass filter.}
   \label{fig:solp_filter}
\end{subfigure}\\
\begin{subfigure}[b]{.45\textwidth}
  \centering
  \includegraphics[width=\textwidth]{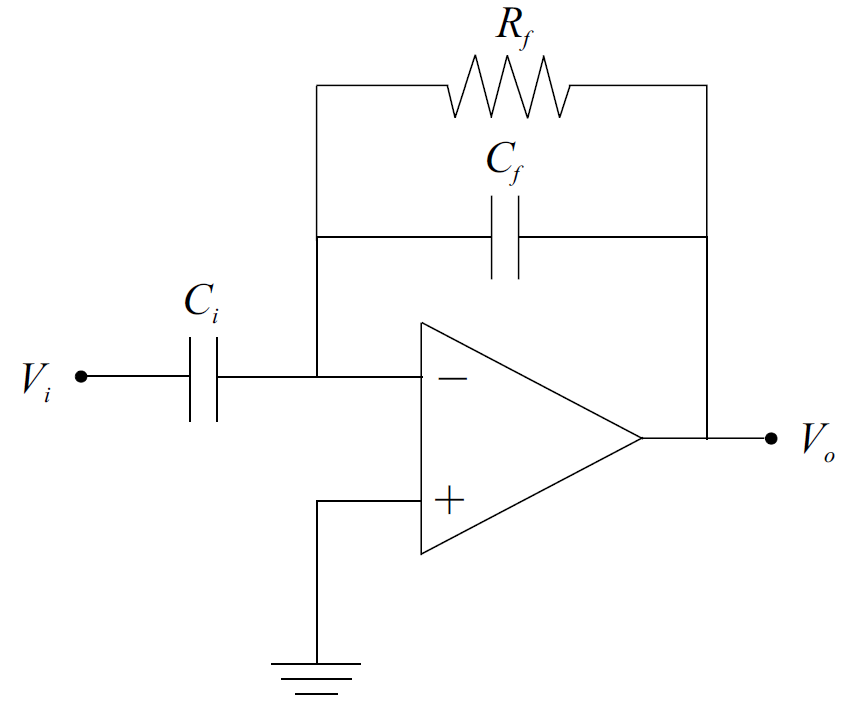}
   \caption{first order high pass filter.}
   \label{fig:fohp_filter}
\end{subfigure}
\begin{subfigure}[b]{.45\textwidth}
  \centering
  \includegraphics[width=\textwidth]{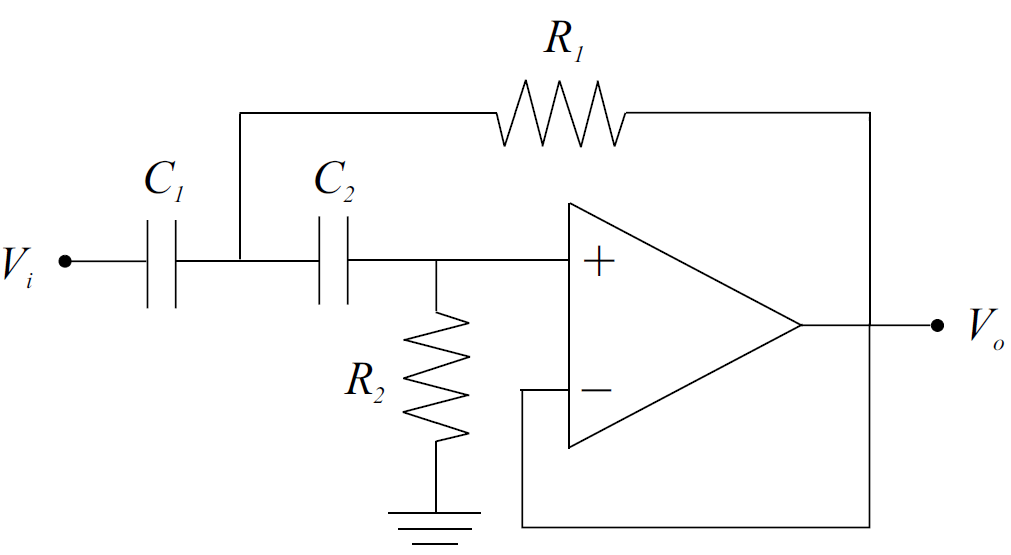}
   \caption{second order high pass filter.}
   \label{fig:sohp_filter}
\end{subfigure}
\caption{First and second order operational amplified based filters.}
\label{fig:op_amp_filters}
\end{figure}

\begin{figure}
\centering
  \centering
  \includegraphics[width=15cm]{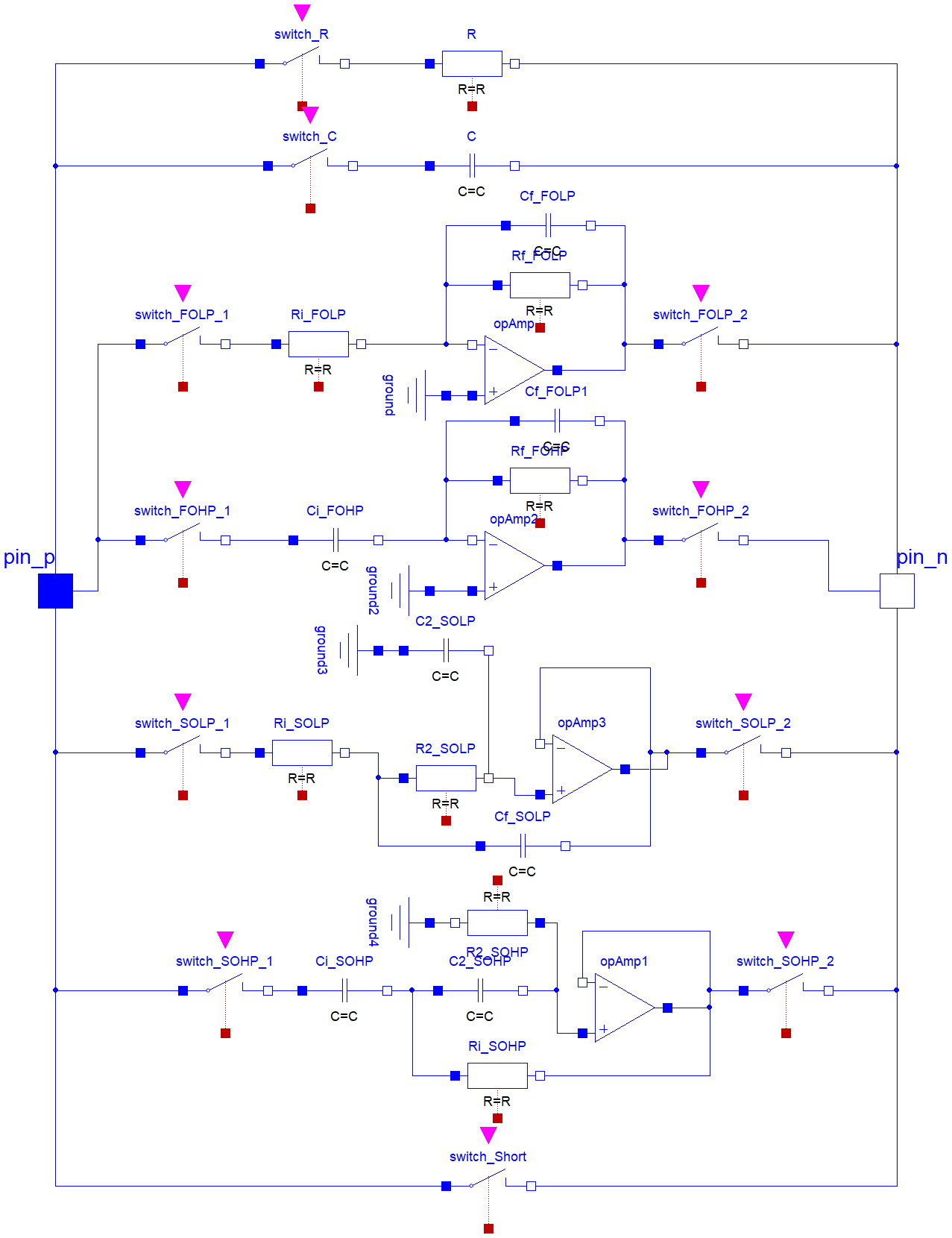}
   \caption{Universal component including low/high, first/second order OpAmp-based filters.}
   \label{fig:universal_component_op_amps_filters}
\end{figure}

In the first approach, we consider basic configurations of OpAmps such as inverting and non-inverting configurations (see Figure \ref{fig:op_amp_configurations}). These configurations result in circuits with only two connectors. To test our design algorithms, we defined a new universal component that removes the inductor and adds inverting and non-inverting OpAmp configurations. We added two switches controlled by a single parameter to indicate the presence or absence of such configurations in the universal component. This parameter is included in the $L_1$ regularization term.
The new modified universal component that includes active components is shown in Figure \ref{fig:universal_component_op_amps}.
In the second approach, we consider even more intricate OpAmp-based atomic components as part of the universal component. For instance, we can incorporate configurations that specify first- or second-order filters, which, when combined in series or parallel, result in higher-order transfer functions. Figure \ref{fig:op_amp_filters} illustrates implementations of such filters using OpAmps. A universal component can be defined by connecting these filters in parallel and equipping each filter with switches for activation.
Figure \ref{fig:universal_component_op_amps_filters} depicts the new universal component that includes low/high, firts/second order OpAmp-based filters.

\subsection{Differential programming for gradient-based optimization}
Both algorithms introduced in the previous sections use gradient-free optimization to search for the component parameters. The advantage of such algorithms is that they work directly with computational representations of the design model (i.e., FMUs). The disadvantage is that they become slower as the number optimization variables increases. An alternative to gradient-free algorithms is gradient-based algorithms, and the optimization problem would translate into a nonlinear program with dynamical constraints. When dealing with design models represented as ODEs, we can map the design optimization problem into a framework that supports AD (e.g., Pytorch \cite{paszke2017automatic} or Jax \cite{jax2018github}), and solve the problem using gradient descent algorithms. Such platforms are endowed with ODE solvers that support AD \cite{NEURIPS2018_69386f6b}.  To formulate the problem in frameworks such as Pytorch or Jax, we first need to extract the equations from the Modelica model of the design. One approach is to generate an XML representation for the DAE using the {\tt dumpXMLDAE} function of the OpenModelica \cite{openmodelica.org:doc:system,openmodelica.org:doc:usersguide} scripting language. Alternatively, we can process the flattened Modelica using a Python Modelica parser such as {\tt modparc} \cite{modparc}. The extracted equations are converted into symbolic objects such as {\tt Sympy} \cite{10.7717/peerj-cs.103} objects, and mapped into deep-learning platform objects that support automatic differentiation. This process leads to a constrained optimization problem that in the case of the continuous relaxation approach is given by:
\begin{eqnarray}
  \label{eq:11212028}
  \min_{\boldsymbol{x}, \boldsymbol{p}, \boldsymbol{s}} & & \mathcal{C}(\hat{\boldsymbol{y}}_{0:T}(\boldsymbol{p}, \boldsymbol{s}),\boldsymbol{y}_{0:T})+\lambda \|\boldsymbol{s}\|_1  \\
    \label{eq:11212029}
  \textmd{subject to: } & & \dot{\boldsymbol{x}} = f(\boldsymbol{x}, \boldsymbol{z};\boldsymbol{p}, \boldsymbol{s}), \\
     \label{eq:11212029_}
   & &  g(\boldsymbol{x}, \boldsymbol{z};\boldsymbol{p}, \boldsymbol{s}) = 0,\\
    \label{eq:11212030}
  & & \hat{\boldsymbol{y}} = h(\boldsymbol{x},  \boldsymbol{z};\boldsymbol{p}, \boldsymbol{s}),
\end{eqnarray}
where (\ref{eq:11212029})-(\ref{eq:11212029_}) is the DAE in the semi-explicit form, representing the dynamics of the design model, and $h(\boldsymbol{x};\boldsymbol{p}, \boldsymbol{s})$ is the sensing model.

To solve (\ref{eq:11212028}), we can convert (\ref{eq:11212029}) into a set of equality constraints using direct collocation methods \cite{doi:10.2514/3.20223,doi:10.2514/3.21662}, or we can use local (e.g., Chebyshev polynomial expansions \cite{boyd01}) of global (e.g., neural networks) parameterizations the state solution and solve for the representation parameters (e.g., weights and biases of the neural network). For example, if we use neural networks to represent the state $\boldsymbol{x}(t)=NN_x(t;\beta_x)$ and the algebraic variables $\boldsymbol{z}(t)=NN_z(t;\beta_z)$, the optimization problem (\ref{eq:11212028}) will be solved in terms of the parameters $\beta_x$, $\beta_z$, $\boldsymbol{p}$, $\boldsymbol{s}$. In addition, automatic differentiation can be used to evaluate the time derivative of the state. Our attempts to use a differentiable programming paradigm to solve design problems were met with mixed results. In the case where the model is represented as an ODE, we obtained good results. For example, in \cite{9147287} we showed how to learn control policies for an inverted pendulum using a model predictive control approach solved using Pytorch. When dealing with DAEs though, the gradient-based optimization algorithm, when combined with direct collocation methods to approximate time derivatives, tend to converge slowly. In addition, the parameterized DAE solution does not always check against the DAE simulation executed with the optimal component and switch parameters. Unfortunately, we cannot always guarantee that the design model can be represented as an ODE, especially since the model is repeatedly reconstructed and simplified. Thus, we prefer to use a direct method (i.e., Powell algorithm), instead a gradient-based approach. Ideally, we would like to have a sensitivity analysis method embedded in the DAE solvers, so that we can access the Jacobian of the state with respect to the model parameters. Such a method is present for instance in the SUNDIALS software family, introduced in \cite{gardner2022sundials,hindmarsh2005sundials}, with DAE solvers such as CVODES and IDAS that include both direct and adjoint-based approaches to compute sensitivities. Currently though, deep learning platform do not offer such a functionality, except for the case where the DAE can be transformed into an ODE. Moreover, even when dealing with ODE, gradient-descent algorithm that include solvers supporting automatic differentiation tend to slow down as the number of optimization parameters increases. We addressed this challenge in \cite{https://doi.org/10.48550/arxiv.2208.12834}, where we showed that block coordinate descent algorithm in combination with direct collocation method speed up training by several order of magnitude. We are currently working on extending this approach to DAE models.

\section{Results}
\label{sec:results}
In this section we present design results based on the proposed approaches for various design examples.
\subsection{Cauer analog low pass filter via continuous relaxation with passive components}
\begin{figure}[!htp]
\centering
\includegraphics[width=15cm]{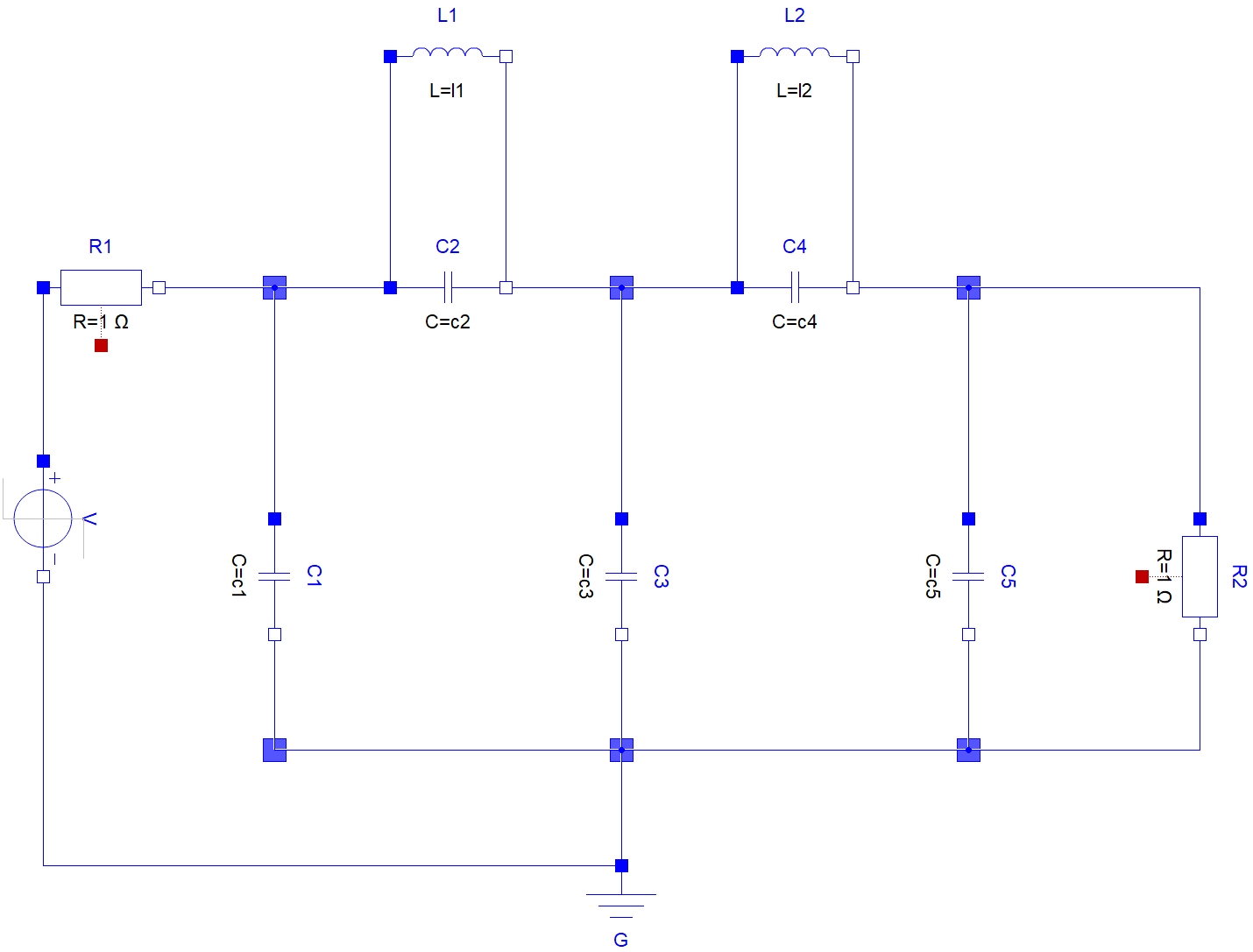}
\caption{Modelica model of the Cauer analog, low pass filter of the fifth order.}
\label{fig:cauer}
\end{figure}

\begin{figure}[!htp]
\centering
\includegraphics[width=15cm]{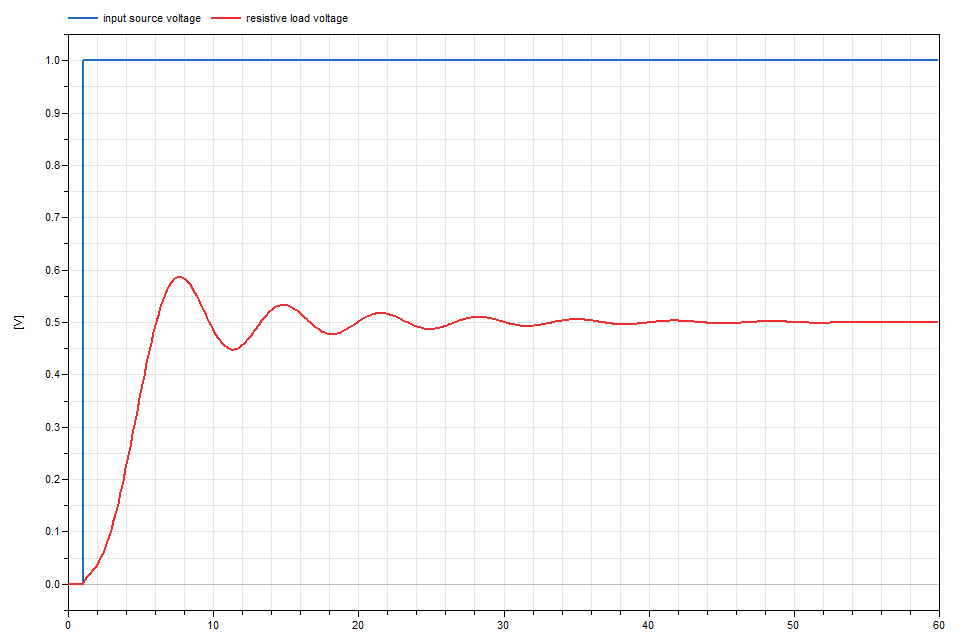}
\caption{Cauer low pass analog filter: input source voltage vs. resistive load voltage.}
\label{fig:cauer_plots}
\end{figure}

\begin{figure}[!htp]
\centering
\includegraphics[width=15cm]{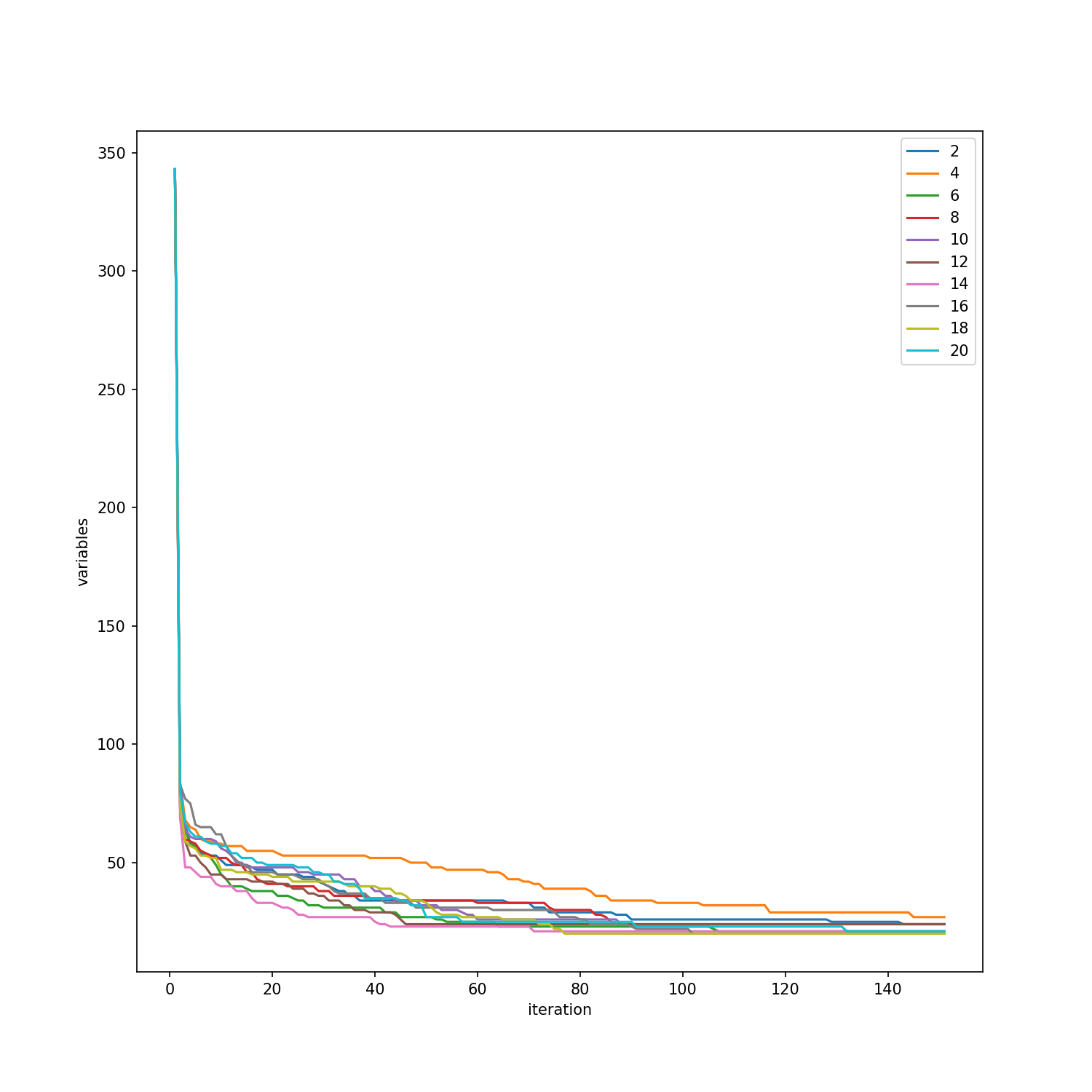}
\caption{Design of the Cauer analog filter with passive components using Algorithm \ref{alg:1}: number of variables decay for even processes.}
\label{fig:num_vars}
\end{figure}

\begin{figure}[!htp]
\centering
\includegraphics[width=15cm]{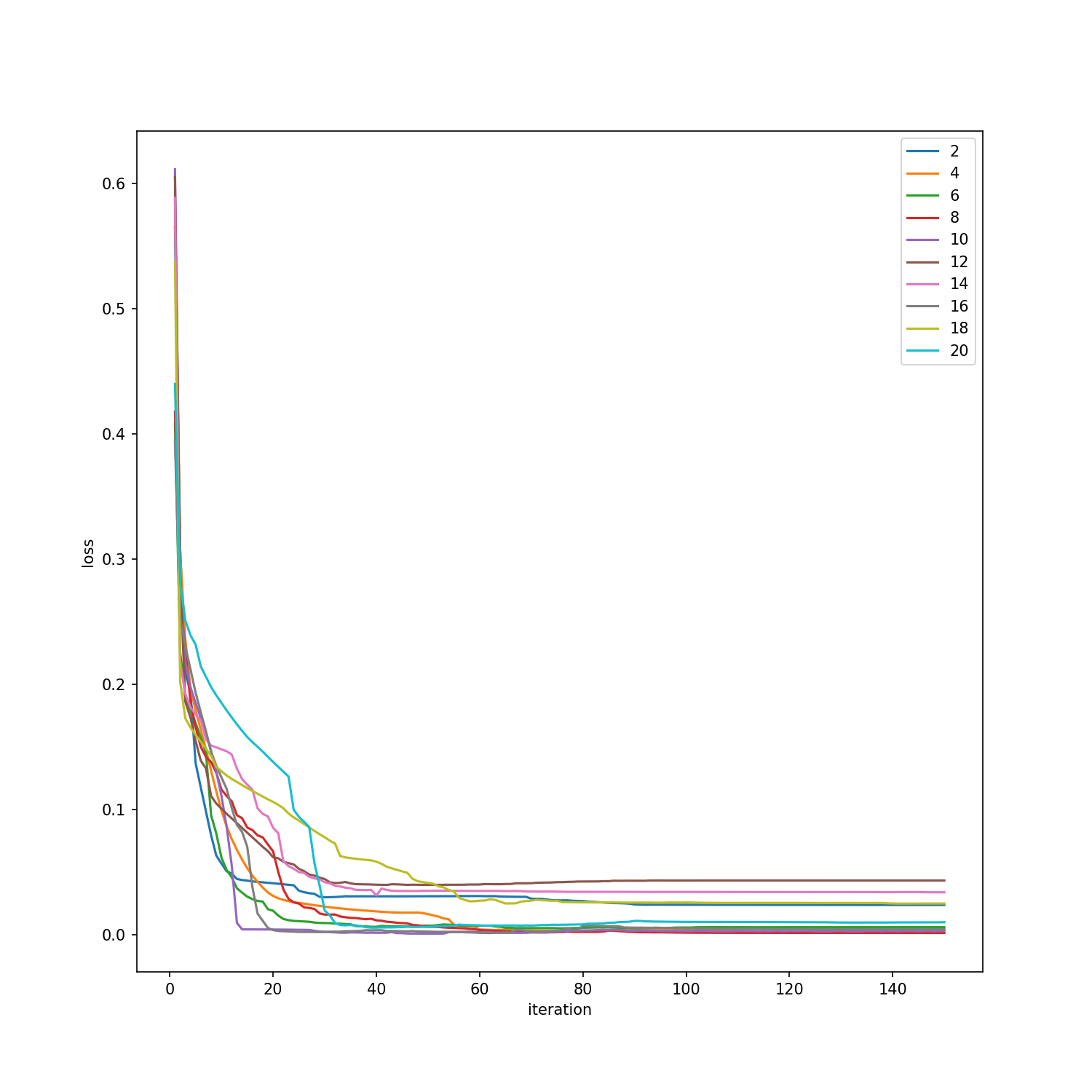}
\caption{Design of the Cauer analog filter with passive components using Algorithm \ref{alg:1}: loss decay for even processes.}
\label{fig:loss}
\end{figure}

\begin{figure}[!htp]
\centering
\includegraphics[width=15cm]{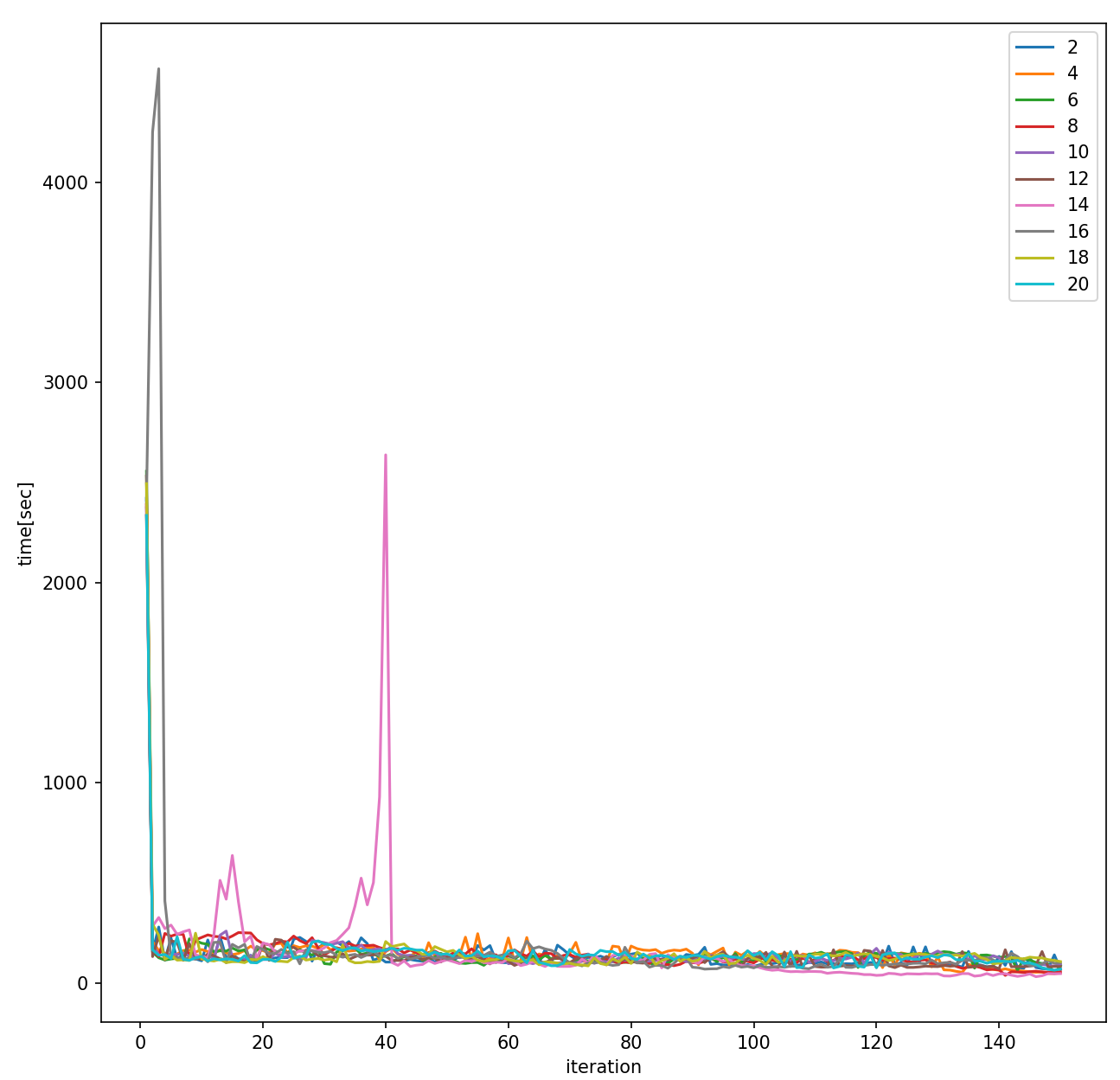}
\caption{Design of the Cauer analog filter with passive components using Algorithm \ref{alg:1}: outer iteration time for even processes.}
\label{fig:iteration time}
\end{figure}

\begin{figure}[!htp]
\centering
\includegraphics[width=15cm]{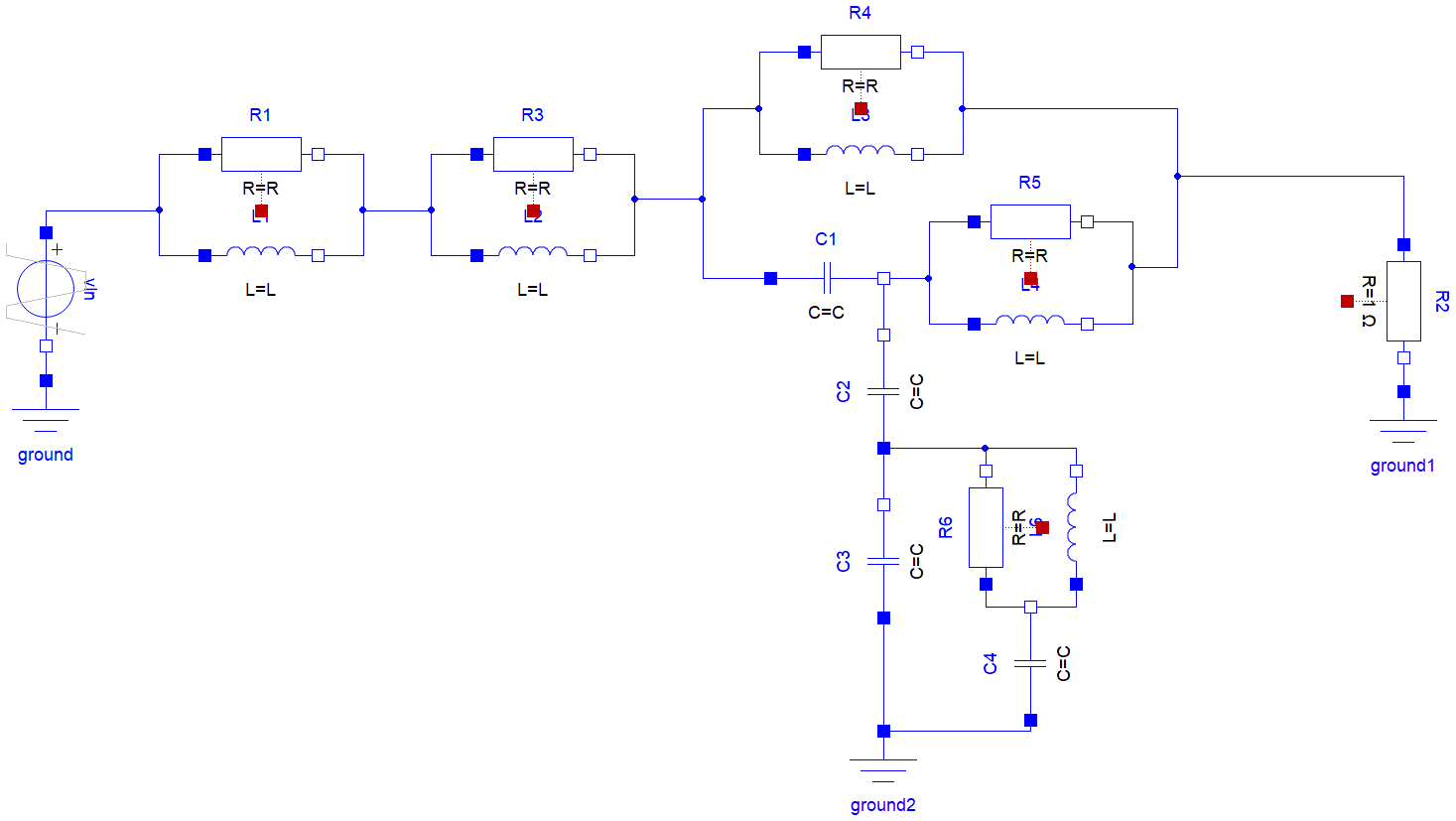}
\caption{Design solution for the Cauer analog low pass filter based on passive components generated by Algorithm \ref{alg:1}.}
\label{fig:Cauer_passive_filter_design_solution}
\end{figure}

Our goal is to design a filter whose output from a step response matches the output of the Cauer analog low pass filter of the fifth order (see Figure \ref{fig:cauer}). The input voltage versus the load voltage plot is shown in Figure \ref{fig:cauer_plots}.
To improve the likelihood to find a design solution, we start with a dense initial topology expressed as a 5x6 grid. The number of optimization variables corresponding to this initial topology is 343, including component parameters and switch values. The dense initial topology is likely to ensure the existence of several local minima that are close to satisfy the design requirements.
\begin{figure}[!htp]
\centering
\includegraphics[width=15cm]{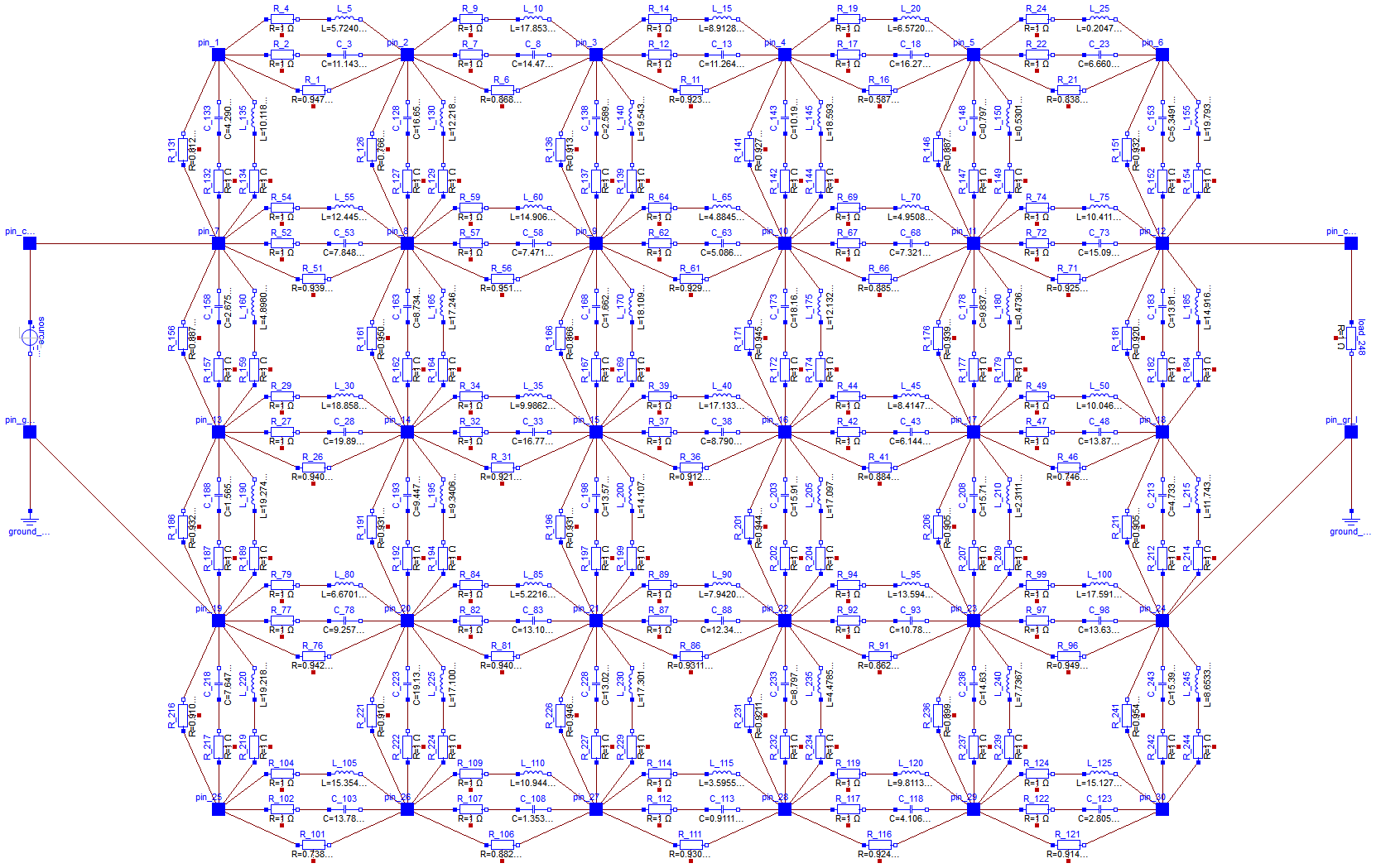}
\caption{Flattened 5x6 grid topology as initial design model.}
\label{fig:flatten_5x6_grid}
\end{figure}
To explore multiple of such local minima, we leverage parallel executions of design optimization processes, where each process starts with random initial component parameters, and initial weight for the $L_1$ cost, and where all switches are initialized to 0.5. We run 20 parallel processes that explore various design solutions. The design optimization algorithm was implemented in Python, and the evaluation of the design loss function was done via FMU-based simulations using the {\tt fmypi} Python package. We refer to each optimization corresponding to an instance of the $L_1$ loss weight as \emph{outer iteration}. An outer iteration was implemented using the gradient free Powell algorithm, where we limit the execution of the algorithm to 150 (inner) iterations. The limited number of iteration affects only the early outer iterations, since 150 iterations may not be sufficient to converge to a local minima. However, since we use a sequence of outer iterations, where each such outer iteration uses the previous optimization variables as initial values, in practice we do converge to a design that satisfies requirements. More importantly, each outer iteration reduces the time complexity since, after each outer iteration we eliminate redundant components whose switch values are approximately zero. This phenomenon can be clearly seen in Figure \ref{fig:num_vars}, where we depict the number of optimization variables as a function the outer iterations, for the even processes. The number of variables drops from 343 at the first iteration to values in the twenties or smaller, at the last iteration. Remarkably, after the first iteration that uses no $L_1$ regularization term, all processes eliminate more than 250 optimization variables as a result of switches being set to zero.
The loss function $\mathcal{C}(\hat{\boldsymbol{y}}_{0:T}(\boldsymbol{p}, \boldsymbol{s}),\boldsymbol{y}_{0:T})$ decay is shown in Figure \ref{fig:loss}. We note that not all processes converge to a very small cost value. This is not surprising since we use random initial conditions and $L_1$ weight loss values. In Figure \ref{fig:iteration time} we show the time per each outer iteration. This time is determined by three factors: the number of iteration of the Powell algorithm, the number of optimization variables and the FMU simulation time. Not unexpectedly, the most expensive outer iteration is the first one, that corresponds to 343 optimization variables.  As the design models become simpler, the outer iteration times reduce to tens of seconds. The graphs show some jumps in time complexity for two processes. Our conjecture is that this due to the numerical instability of the FMU while the optimization algorithms explore values for the component parameters.
An example of a design solution that realizes the behavior of the Cauer analog filter implemented using passive components is shown in Figure \ref{fig:Cauer_passive_filter_design_solution}.

\subsection{Cauer analog low pass filter via pseudo-random search with passive components}
\begin{figure}[!htp]
\centering
\includegraphics[width=15cm]{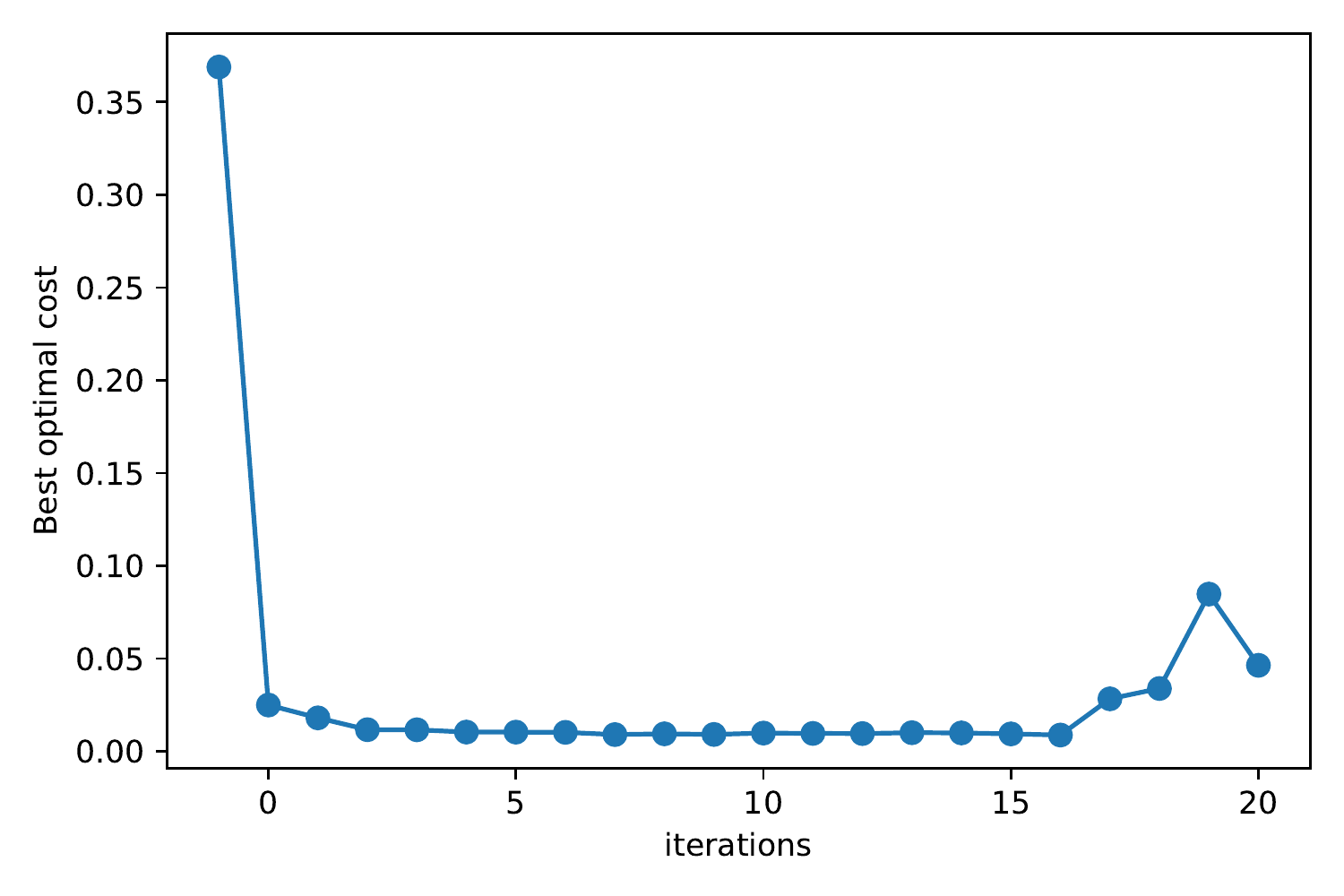}
\caption{Design of the Cauer analog filter with passive components using Algorithm \ref{alg:2}: the minimum cost over optimized topologies (line 8 of Algorithm \ref{alg:2}).}
\label{fig:random_search_optimal_loss}
\end{figure}

\begin{figure}[!htp]
\centering
\includegraphics[width=15cm]{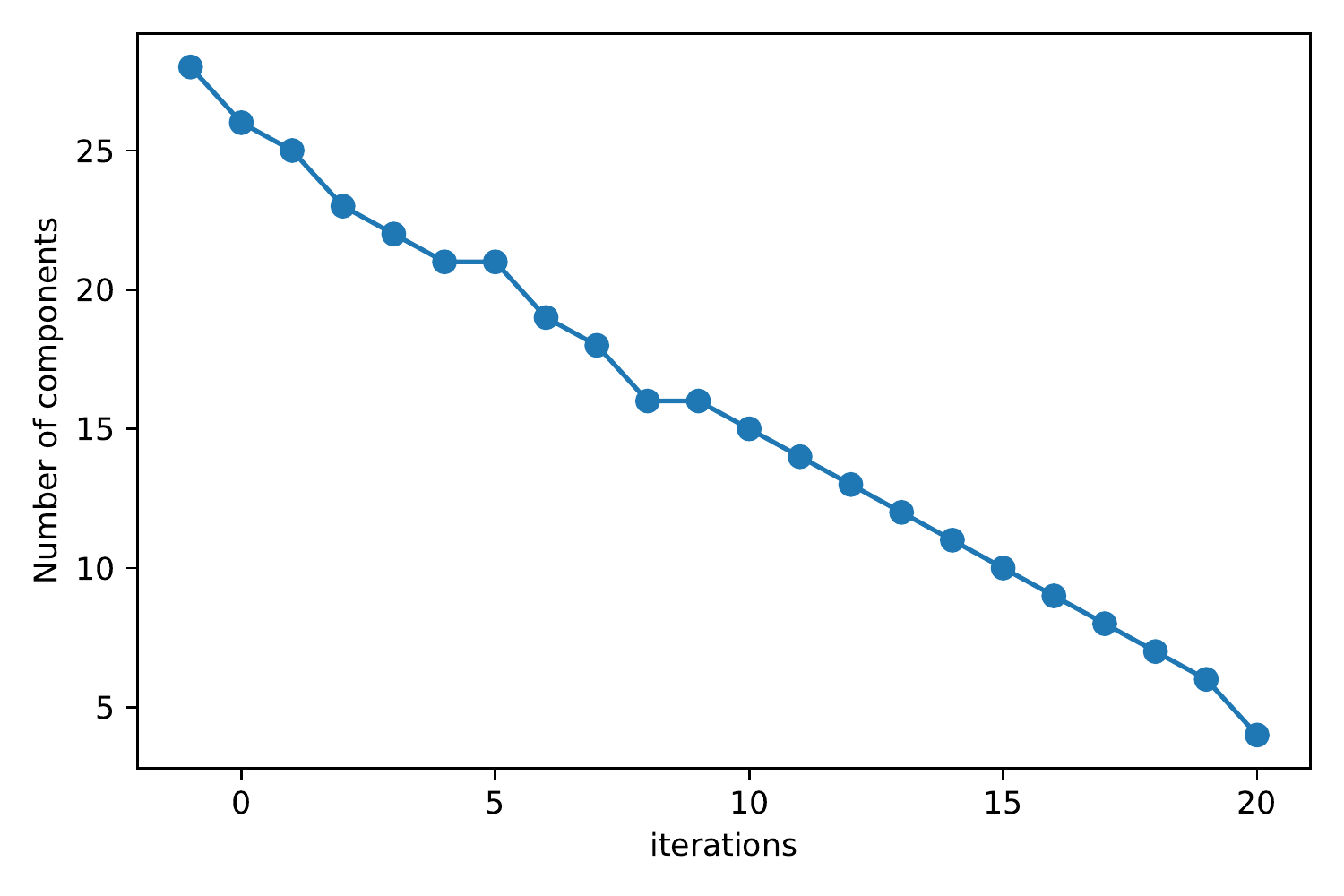}
\caption{Design of the Cauer analog filter with passive components using Algorithm \ref{alg:2}: the number of components corresponding to the topology with the smallest optimal cost.}
\label{fig:random_search_num_components}
\end{figure}

\begin{figure}[!htp]
\centering
\includegraphics[width=15cm]{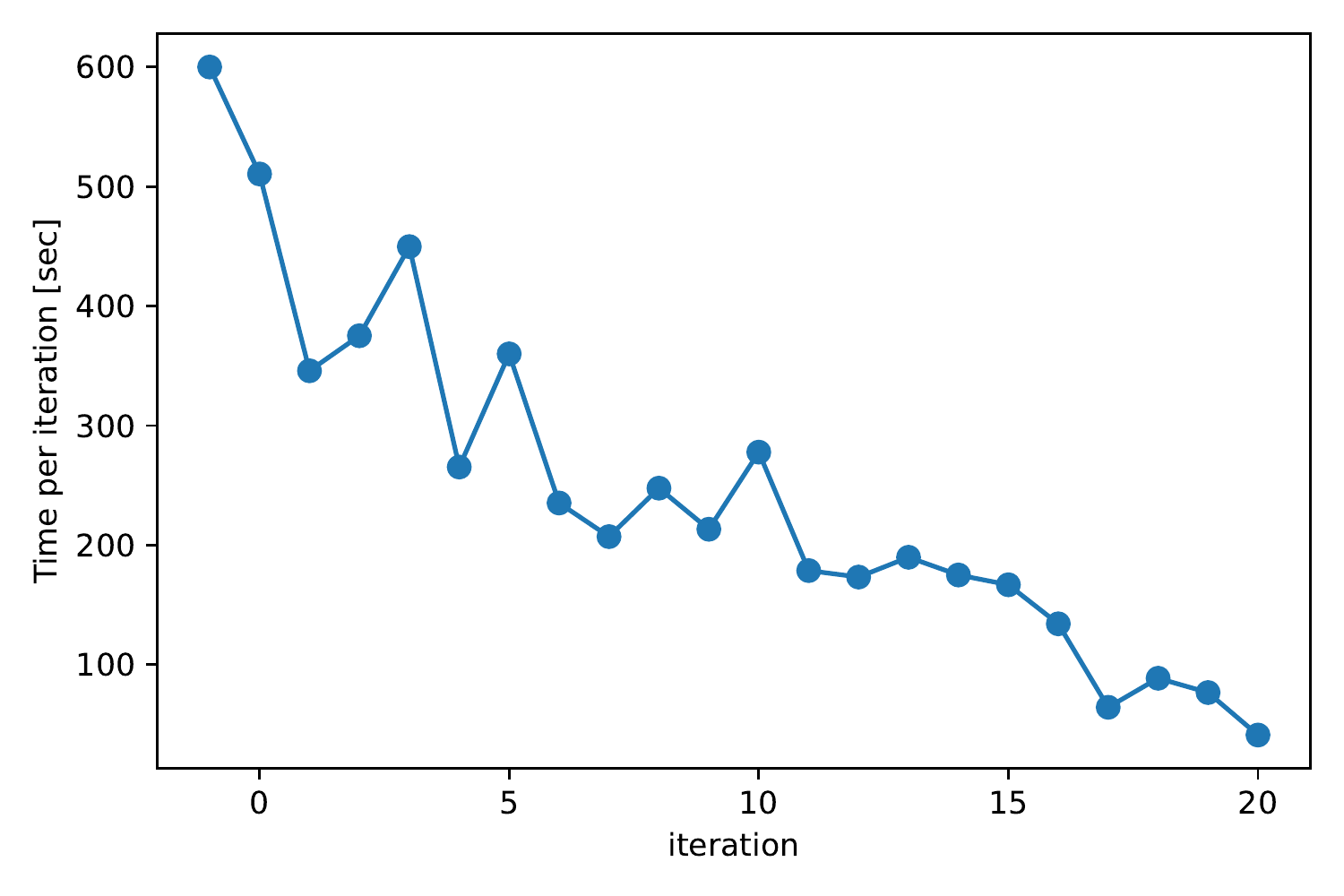}
\caption{Design of the Cauer analog filter with passive components using Algorithm \ref{alg:2}: the duration of each outer iteration (lines 4-24 of Algorithm \ref{alg:2}).}
\label{fig:random_search_iteration_time}
\end{figure}

As in the continuous relaxation case, the objective is to construct an electric circuit whose output matches the output of the Cauer analog low pass filter shown in Figure \ref{fig:cauer}. We use the same design model template given by a 5x6 grid of universal components, as shown in Figure \ref{fig:5x6_grid}. The main difference is that each universal component will be in only one of the five modes, thus the complexity of the randomly generated initial designs will be lower. We start the algorithm with 200k randomly generated topologies based on the 5x6 grid, and simulate them in parallel. The simulation time for all 200k design models took about 5 minutes. At this time the execution of the outer iterations (lines 4-24 of Algorithm \ref{alg:2}) commences, where $n_o=25$, and $C_{th} =0.1$. We solve the parameter optimization problem (line 7 of Algorithm \ref{alg:2}) using Powell algorithm with 500 maximum number of iterations. We leverage parallel process execution to solve $n_o$ optimization problems. Figure \ref{fig:random_search_optimal_loss} shows the best optimized requirements cost for each outer iteration (line 8 of Algorithm \ref{alg:2}), while Figure \ref{fig:random_search_num_components} shows the number of components corresponding to the topologies with the best costs. We note a continuous decrease of the cost, until the number of components becomes too small for the requirements to be satisfied. The number of components just before the increase of the cost is 9, where the number of components in the Cauer circuit (excluding the load resistance) is 8. Figure \ref{fig:random_search_iteration_time} shows the time per each outer iteration that includes both model simulations (line 3 of Algorithm \ref{alg:2}) and parameter optimizations (line 7 of Algorithm \ref{alg:2}). The reduced complexity of the models is emphasized by the decrease of the outer iteration time. Note that, unlike Algorithm \ref{alg:1}, this effect is not obtained via model reconstruction, but only via the reduction in the number of model parameters and  equations. As components are eliminated via {\tt short} or {\tt open} modes, the number of equations decreases.

\begin{figure}[!htp]
\centering
\includegraphics[width=15cm]{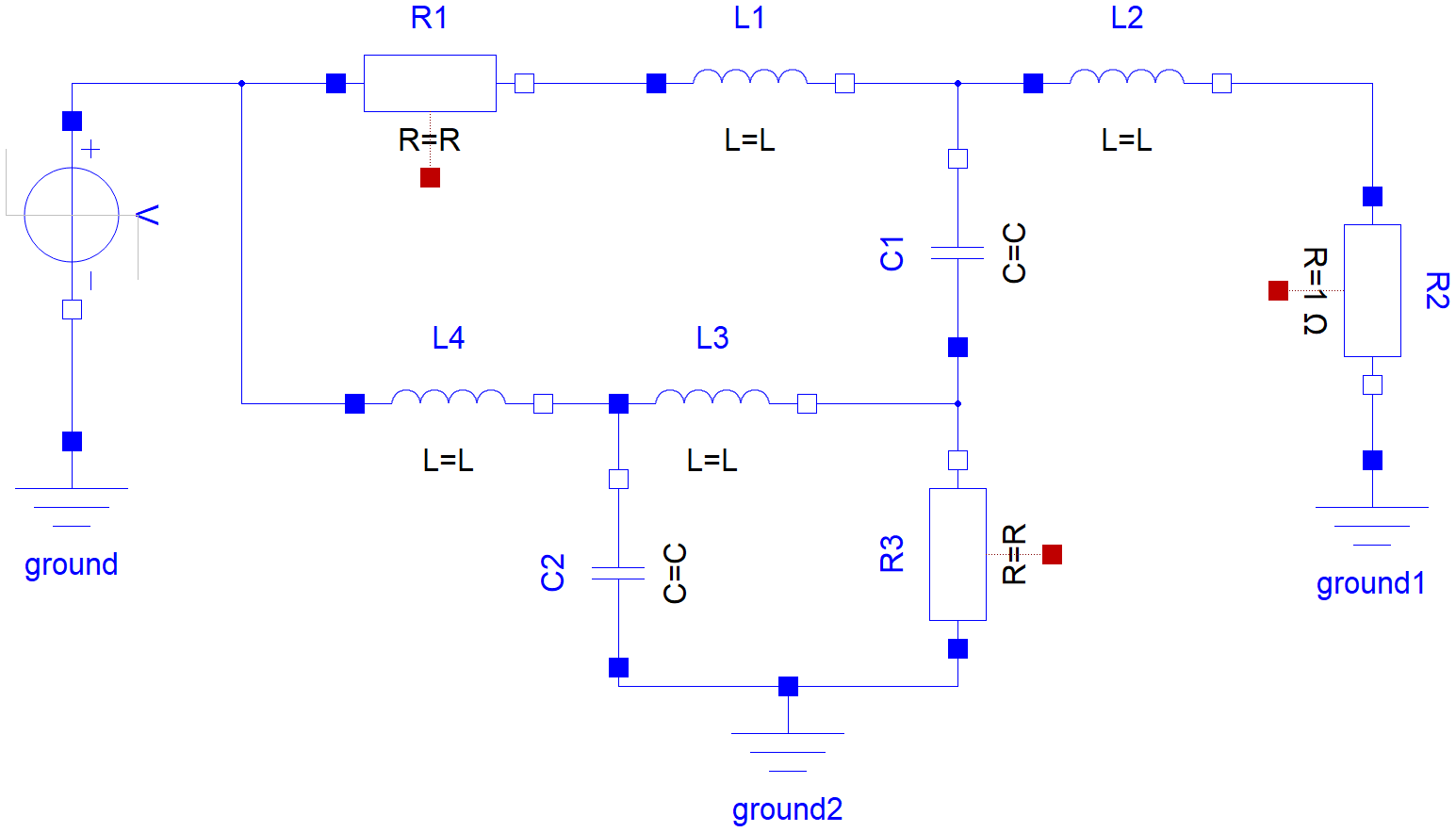}
\caption{Design solution for the Cauer analog low pass filter based on passive components generated by Algorithm \ref{alg:2}}
\label{fig:cauer_design_solution_random_search}
\end{figure}

Figure \ref{fig:cauer_design_solution_random_search} shows the realization of the best trade-off between the number of components and accuracy cost of the Cauer analog filter using passive components.

\subsection{Voltage level shifter design via continuous relaxation with operational amplifier}
\begin{figure}[!htp]
\centering
\includegraphics[width=15cm]{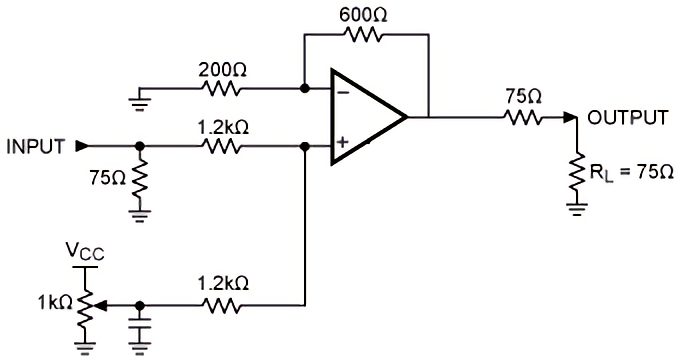}
\caption{Voltage level shifter circuit used to generate the ground truth data in the form of the voltage across the load resistor ($R_L$). }
\label{fig:voltage_shifter}
\end{figure}

We present the results of designing a voltage level shifter(see Figure \ref{fig:voltage_shifter}), using Algorithm \ref{alg:1}. The universal component employed to generate the initial grid topology consists of a resistor, capacitor, and operational amplifier arranged in a non-inverting configuration.. We executed 10 parallel processes of Algorithm \ref{alg:1} for 150 outer iterations, with a limit of 300 inner iterations for the Powell algorithm in each outer iteration.

\begin{figure}[!htp]
\centering
\includegraphics[width=15cm]{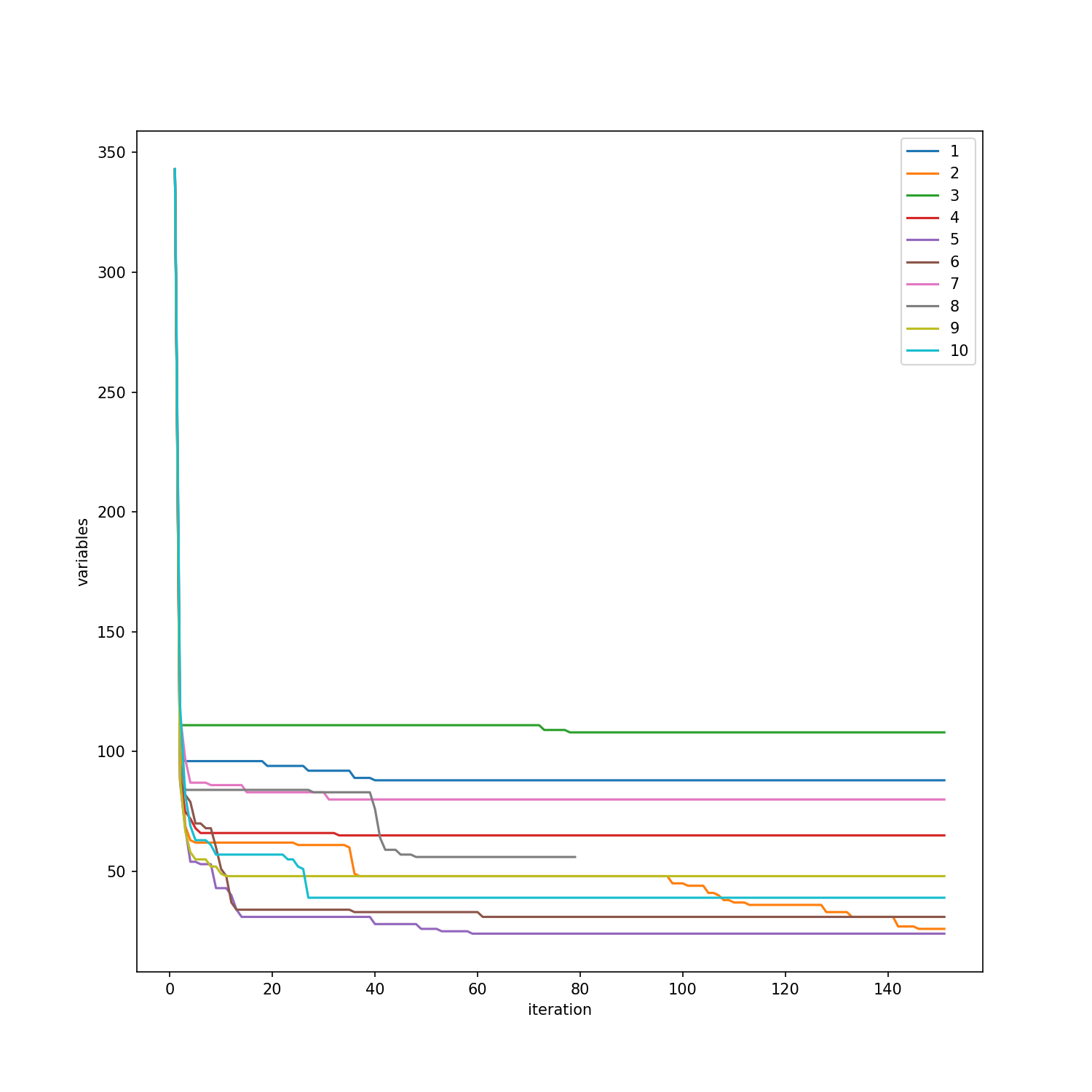}
\caption{Design of the voltage level shifter  with operational amplifiers using Algorithm \ref{alg:1}: number of variables decay for even processes.}
\label{fig:num_vars_level_shifter}
\end{figure}

\begin{figure}[!htp]
\centering
\includegraphics[width=15cm]{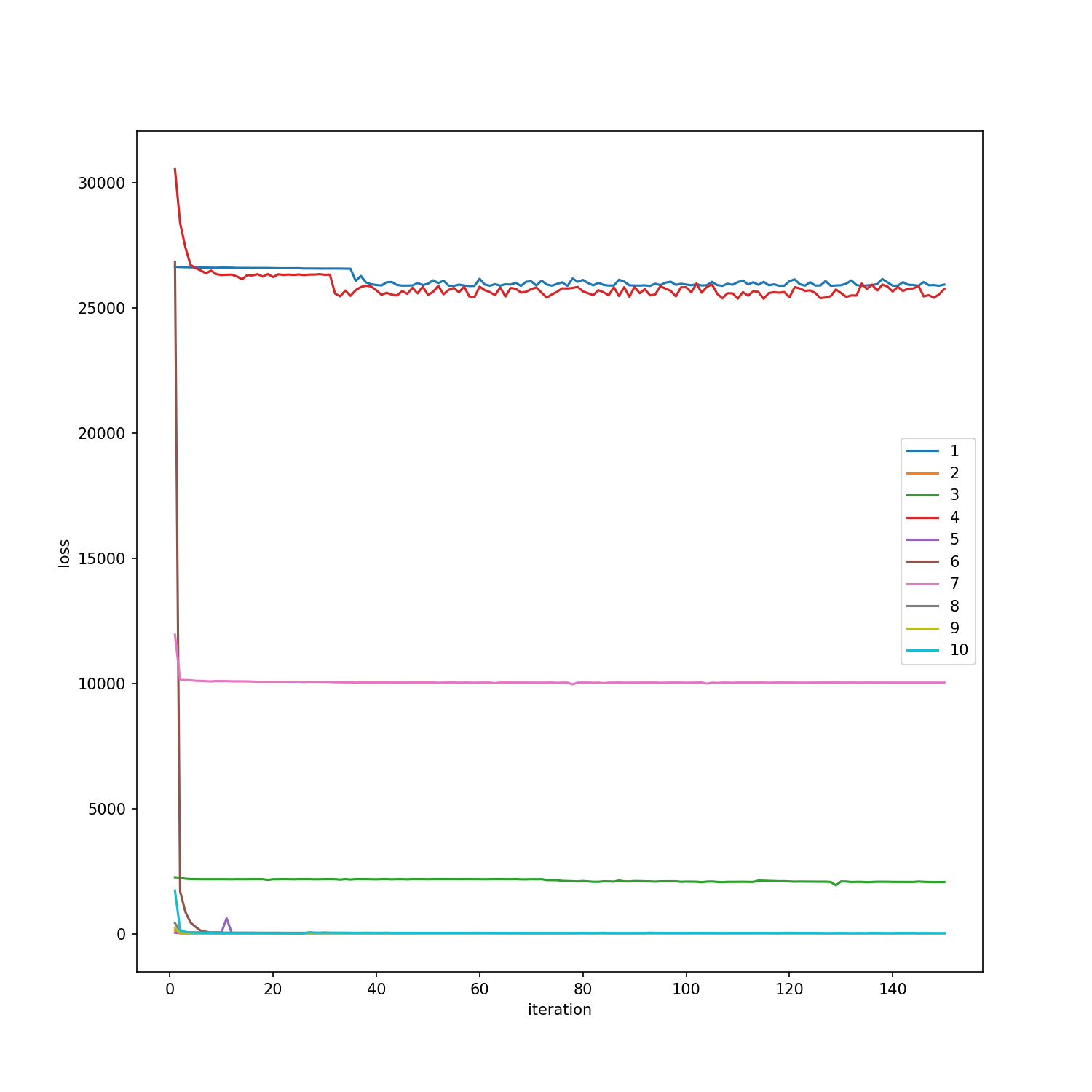}
\caption{Design of the voltage level shifter  with operational amplifiers using Algorithm \ref{alg:1}: loss decay for even processes.}
\label{fig:loss_level_shifter}
\end{figure}

\begin{figure}[!htp]
\centering
\includegraphics[width=15cm]{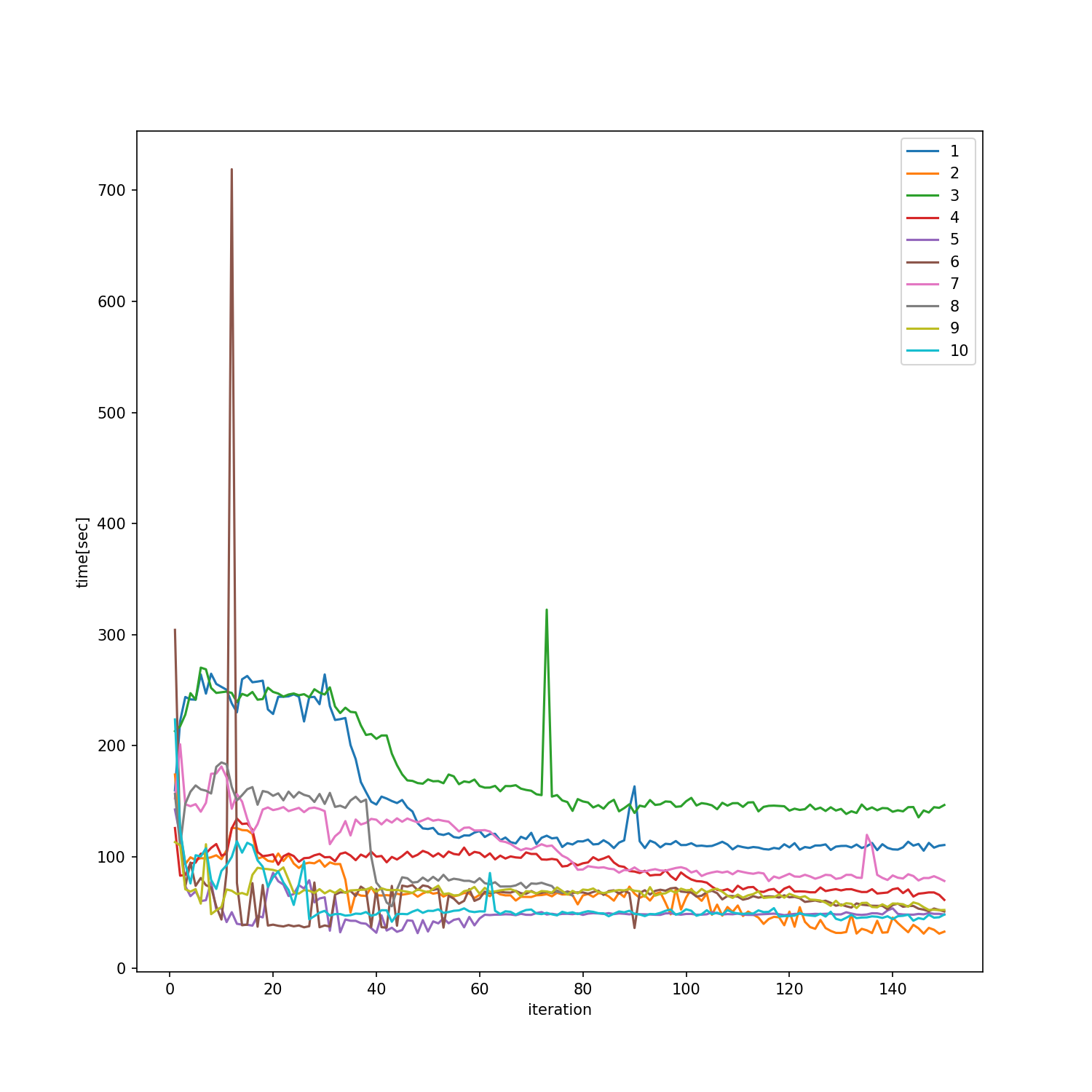}
\caption{Design of the voltage level shifter  with operational amplifiers using Algorithm \ref{alg:1}: outer iteration time for even processes.}
\label{fig:iteration time_level_shifter}
\end{figure}

\begin{figure}[!htp]
\centering
\includegraphics[width=15cm]{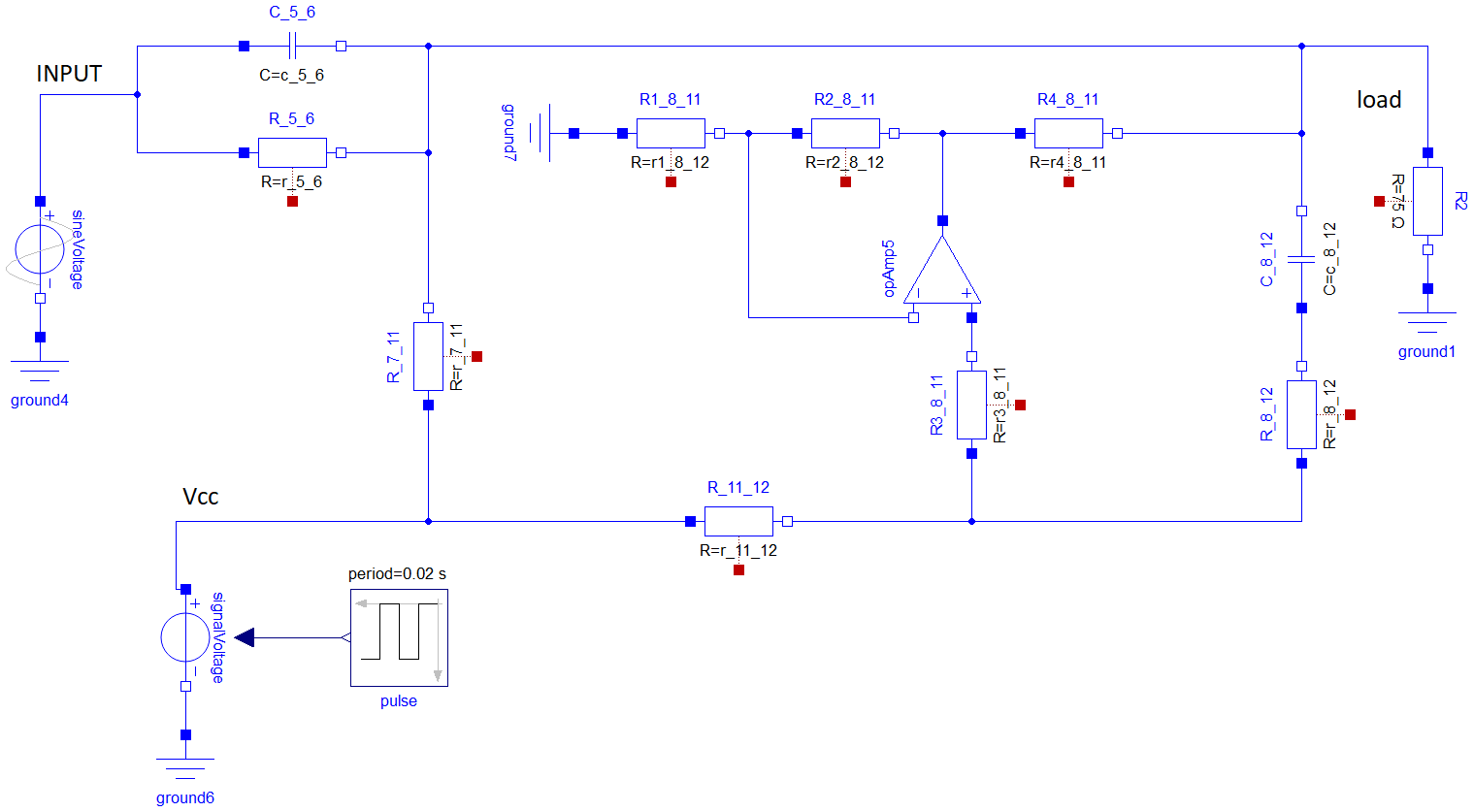}
\caption{Design of the voltage level shifter  with operational amplifiers using Algorithm \ref{alg:1}: design solution 1.}
\label{fig:voltage level shifter solution 1}
\end{figure}

\begin{figure}[!htp]
\centering
\includegraphics[width=15cm]{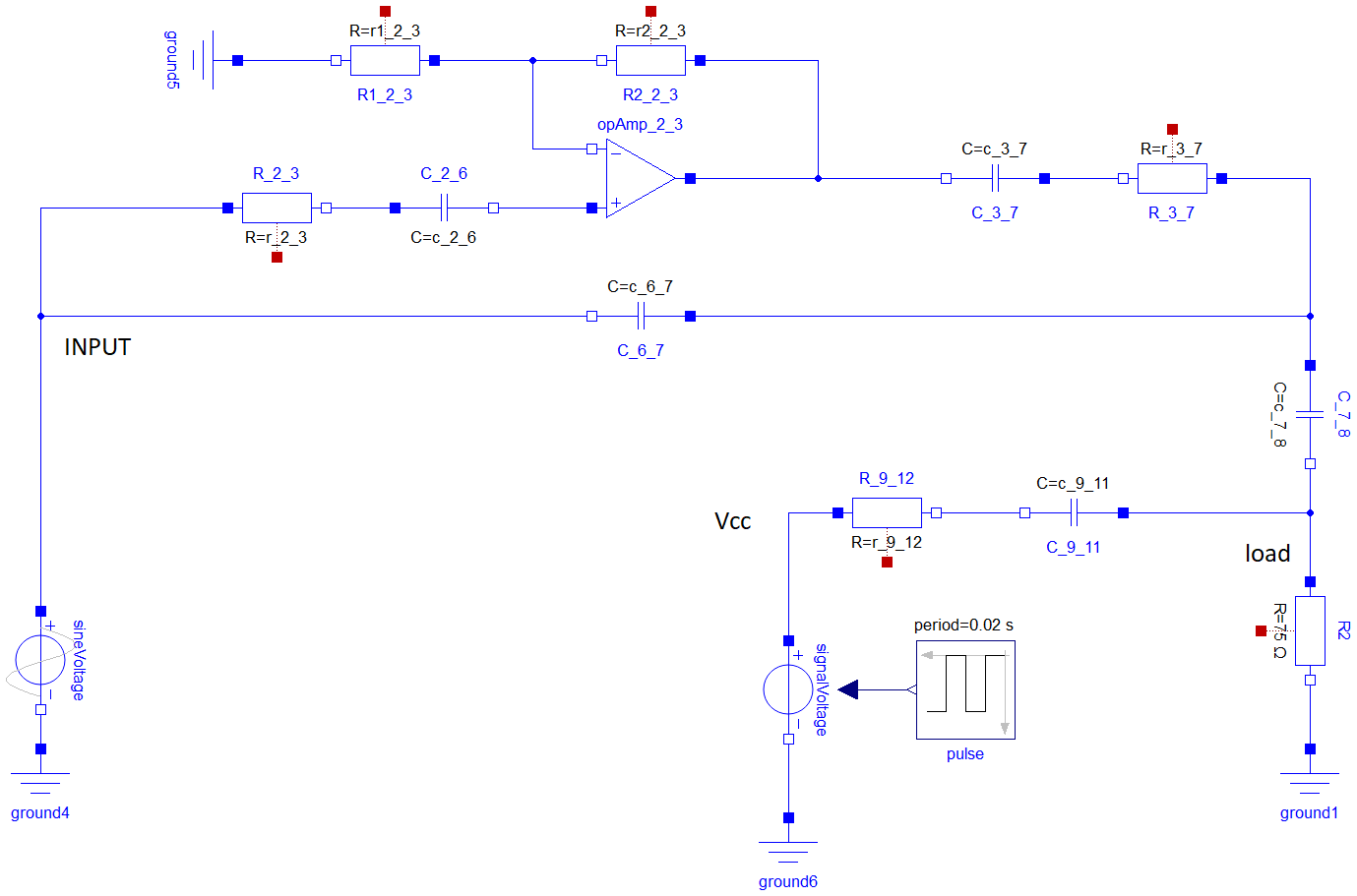}
\caption{Design of the voltage level shifter  with operational amplifiers using Algorithm \ref{alg:1}: design solution 2.}
\label{fig:voltage level shifter solution 2}
\end{figure}

The learning results for the 10 processes are depicted in Figures \ref{fig:num_vars_level_shifter}, \ref{fig:loss_level_shifter} and  \ref{fig:iteration time_level_shifter}, respectively, which illustrate the evolution of the number of optimization variables, the loss function for requirements, and the execution time for the outer iterations.
The results indicate that processes 5 and 6 yield the most favorable outcomes in terms of both the number of optimization variables and loss function values. Based on these two processes, the design solutions for the voltage level shifter are depicted in Figures \ref{fig:voltage level shifter solution 1} and \ref{fig:voltage level shifter solution 2}.
Notably, both solutions have a component count that is similar to that of the original level shifter depicted in Figure \ref{fig:voltage_shifter}, with 10 and 9 components for the two solutions compared to 8 components in the original circuit (not counting the load resistor and the voltage source components). Additionally, both solutions utilize a single OpAmp.

\subsection{Cauer analog low pass filter via continuous relaxation with active filters}
\begin{figure}[!htp]
\centering
\includegraphics[width=15cm]{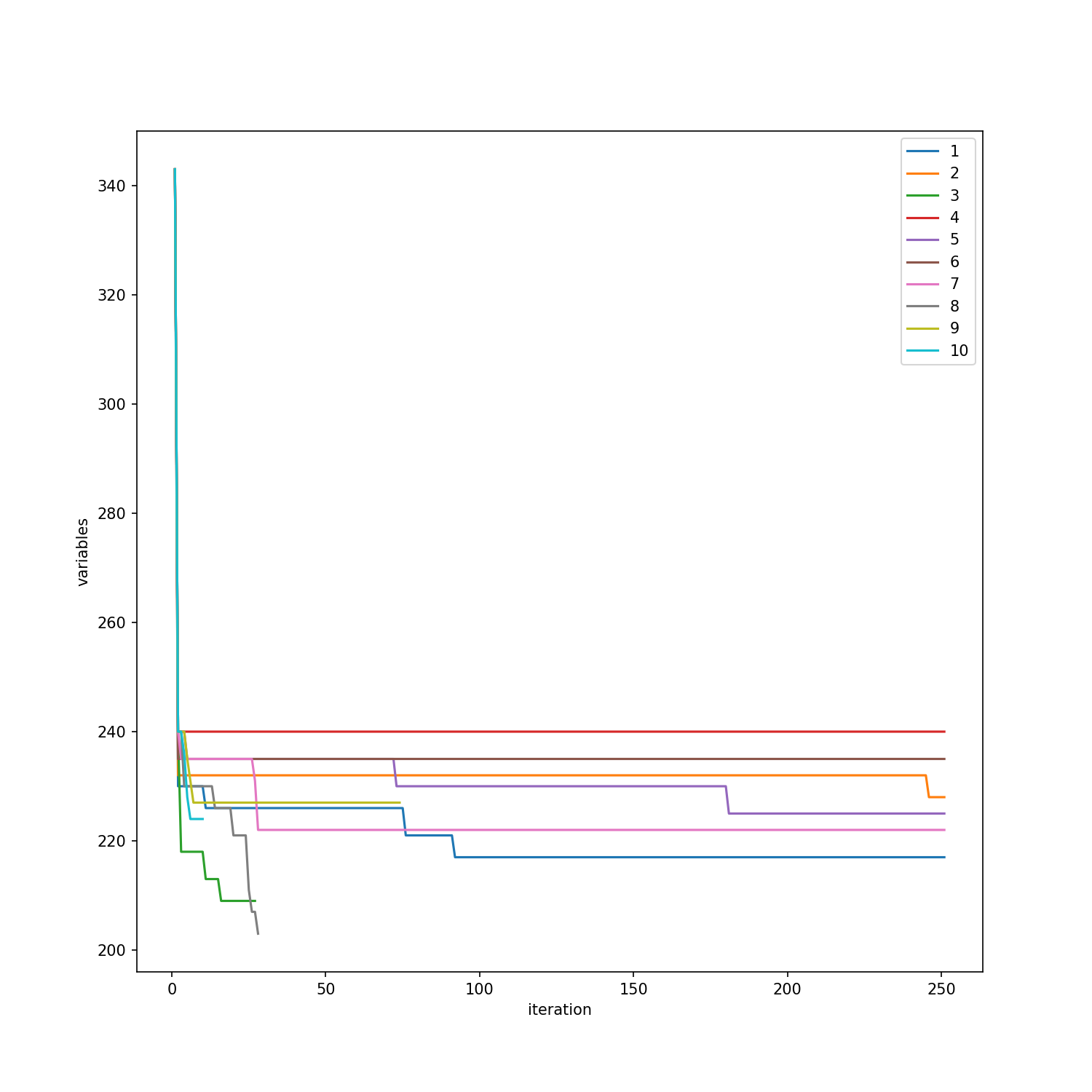}
\caption{Design of the Cauer analog low pass filter with operational amplifier-based first and second order filters, using Algorithm \ref{alg:1}: number of variables decay for even processes.}
\label{fig:num_vars_Cauer_active_filters}
\end{figure}

\begin{figure}[!htp]
\centering
\includegraphics[width=15cm]{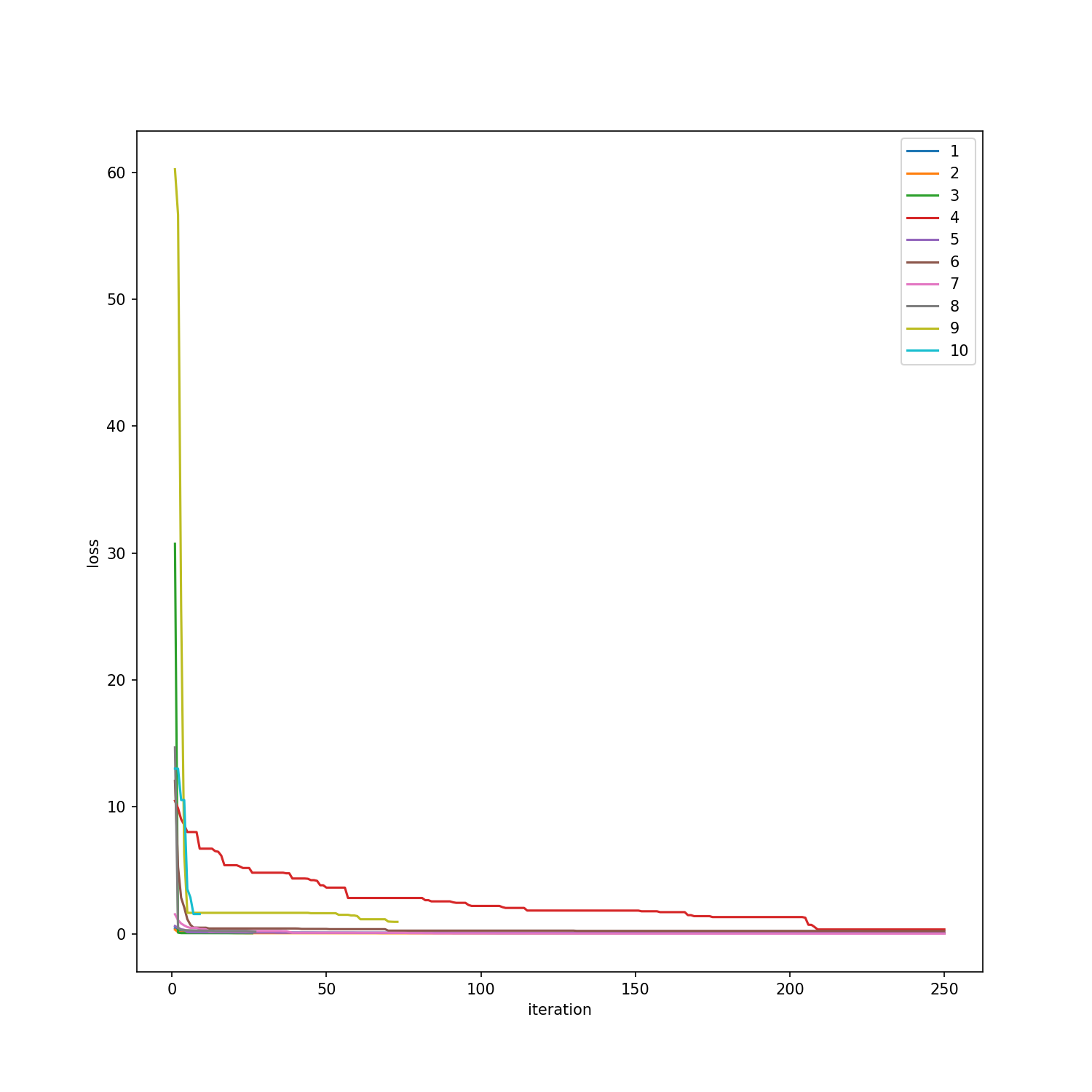}
\caption{Design of the Cauer analog low pass filter with operational amplifier-based first and second order filters, using Algorithm \ref{alg:1}: loss decay for even processes.}
\label{fig:loss_num_vars_Cauer_active_filters}
\end{figure}

\begin{figure}[!htp]
\centering
\includegraphics[width=15cm]{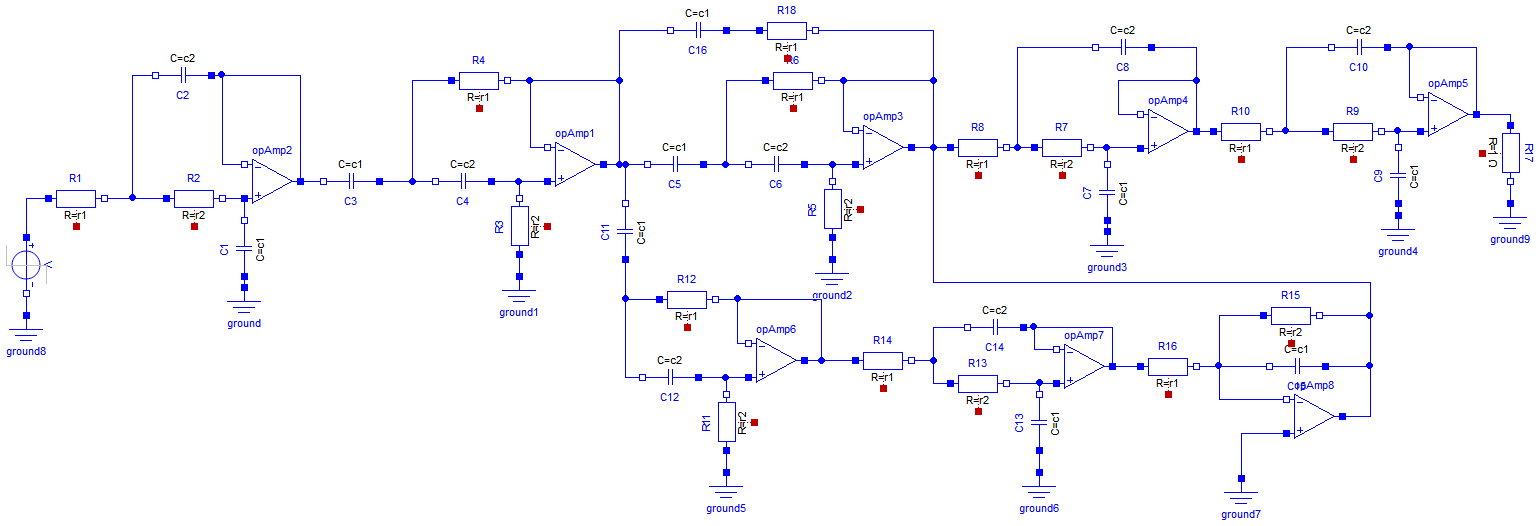}
\caption{Design solution for the Cauer analog low pass filter  with operational amplifiers using Algorithm \ref{alg:1}.}
\label{fig:Cauer_active_filter_design_solution}
\end{figure}

We repeated the design optimization problem for the Cauer low pass filter, where the universal component was defined using a combination of first and second-order low and high-pass filters, implemented using operational amplifiers. We started with a 2x6 grid as the initial topology and ran 10 parallel processes of Algorithm 1 for 250 outer iterations, with a limit of 1000 inner iterations for the Powell algorithm in each iteration. The learning results for the 10 processes are displayed in Figures \ref{fig:num_vars_Cauer_active_filters}, \ref{fig:loss_num_vars_Cauer_active_filters}, and \ref{fig:loss_num_vars_Cauer_active_filters}. Some processes stopped or were slow due to numerical instabilities, but we still generated several design solutions. After a final simplification, we chose one of the solutions and arrived at a circuit shown in Figure \ref{fig:Cauer_active_filter_design_solution} that includes 8 operational amplifiers. The Modelica Standard Library (MSL) has an implementation of the Cauer analog filter that uses only 5 operational amplifiers but also includes 4 negative resistors, where each negative resistor can be implemented using an operational amplifier. Our design solution therefore has a similar number of operational amplifiers as the one in the MSL.
\begin{figure}[!htp]
\centering
\includegraphics[width=15cm]{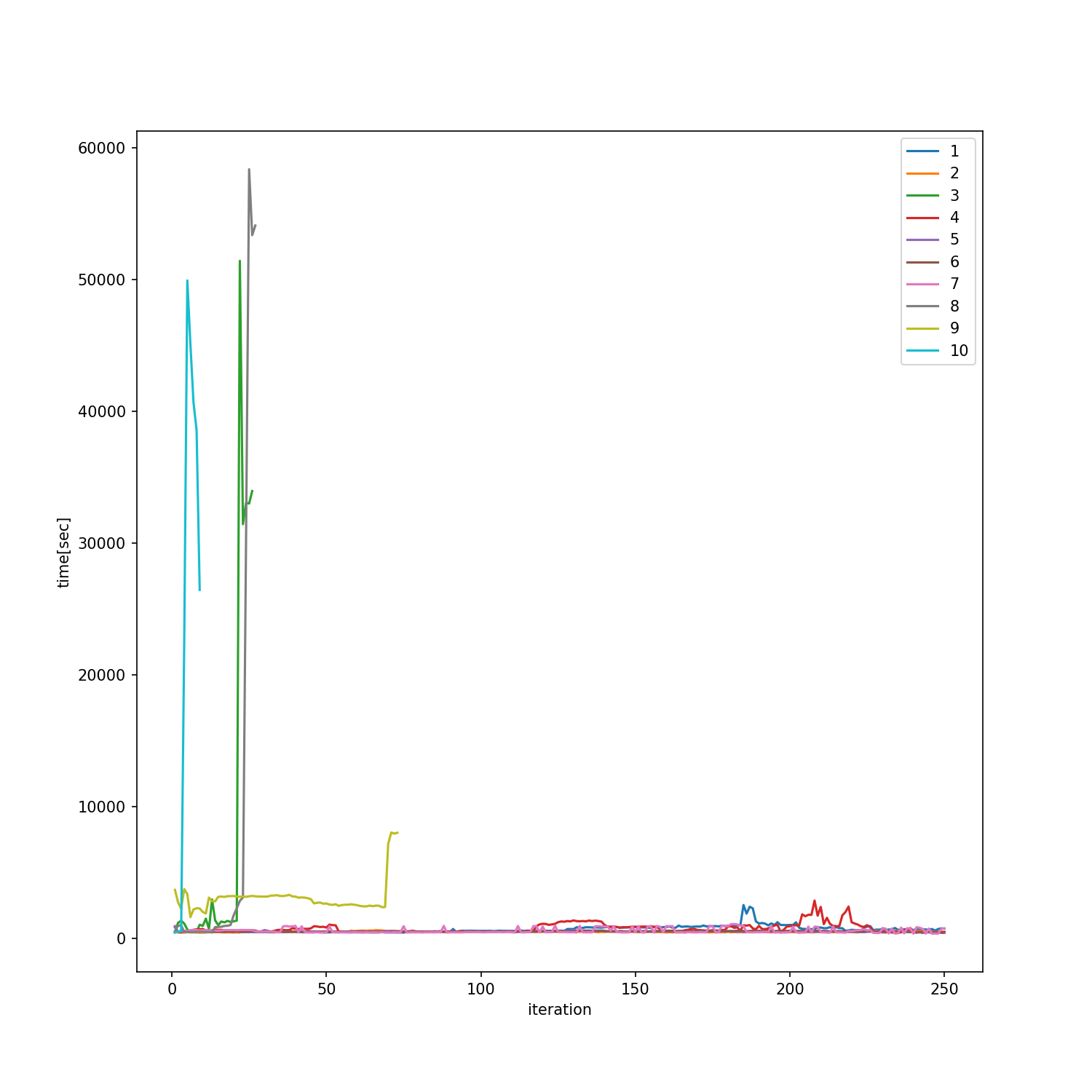}
\caption{Design of the Cauer analog low pass filter with operational amplifier-based first and second order filters, using Algorithm \ref{alg:1}: outer iteration time for even processes.}
\label{fig:iteration time_num_vars_Cauer_active_filters}
\end{figure}

Table \ref{table:1} summarizes the design results of the above examples in comparison with the original circuits that were used to generate the ground truth. When counting the number of resistors and OpAmps in the MSL active implementation of the Cauer filter, we included the number of resistors and OpAmps needed to implement the negative resistors.

\begin{table}[h!]
\centering
\begin{tabular}{|c|c|c|c|c|}
 \hline
 Circuit & Number of resistors & Number of capacitors & Number of inductors & Number of OpAmps \\
 \hline
 Original passive Cauer filter & 1 & 5 & 2 & 0 \\
 \hline
 Designed passive Cauer filter (Alg. 1) & 5 & 3 & 5 & 0 \\
  \hline
 Designed passive Cauer filter (Alg. 2) & 2 & 2 & 4 & 0 \\
   \hline
 Original active Cauer filter  & 19 & 8 & 0 & 9 \\
 \hline
 Designed active Cauer filter (Alg. 1) & 17 & 16 & 0 & 8 \\
 \hline
 Original voltage level shifter  & 7 & 1 & 0 & 1 \\
 \hline
 Designed voltage level shifter (sol. 1)  & 8 & 2 & 0 & 1 \\
 \hline
  Designed voltage level shifter (sol. 2)  & 5 & 4 & 0 & 1 \\
 \hline
\end{tabular}
\caption{Summary of the design results for various examples.}
\label{table:1}
\end{table} 
\subsection{Related work}
The existing work focuses on tools and languages for supporting designs (i.e., Simscape/Matlab, Modelica/OpenModelica), and optimizing components in preset design topologies. Tools supporting design processes are not endowed with automated design features.  It requires the conceptual design as input, where the designer uses their experience to come up with several design solution that are refined through simulations. An approach that combines the Modelica language and model-based system design verification against requirements is discussed in \cite{10555523466162346647}. Other approaches use model templates to create design models that are tested against requirements via simulations \cite{clima_2022}.
In \cite{9867163}, the authors use a component-based approach to multirotor aircraft design. The components were modeled using
a graph-based modeling approach and combined into a system model using automated graph combination tools. In our case, we get compositionality automatically since the Modelica language has a rich semantics that enables composition of physics-based components. Similar to our approach, the authors of \cite{9867163} use a hybrid approach to design optimization. They map a discrete problem to a continuous optimization problem and based on the solution of the latter, they search for a solution in the discrete space. One of the main differences as compared to our setup is that they have the (meta-design) already determined, i.e., the design concept is set. In this scenario, the goal is to find instantiations of components in a fixed design. Our goal is to find the meta-design and component instantiations, jointly. The authors of \cite{9867163} use gradient-based methods to solve the continuous optimization problem, where the gradients are approximated via finite-difference. Such an approach is not scalable when the number of optimization parameters is large. A more effective approach is to use automatic differentiation, since the authors have direct access to the model equations. Typically, design problems attempt to optimize multiple requirements simultaneously. Such an objective leads to multi-objective optimization problems where the goal is to generate a Pareto-optimal front (POF). In \cite{10.1145/3449726.3463194}, the authors show how to use the Tigon optimization library \cite{DURO2021106851} to implement multi-objective evolutionary algorithms that can generate approximations of POF. Our formulation includes multiple objectives as well. However, we decided to follow a lexicographic method and favor one of the objective, and then continuously optimize the second objective (i.e., sparsity objective) as long as the primary objective is not affected. In \cite{Pettersson18090} the authors introduce methods and algorithms for the development of industrial robots. They show how to select the size of the gearboxes and arm lengths using the Complex-RFD \cite{1993631} optimization method. As other design examples, the approach assumes that the topology is already chosen, and the objective is to choose instantiations of components for the fixed topology.
In \cite{10.5555/3437539.3437740} the authors use  graph convolutional neural networks and reinforcement learning (RL) for transistor sizing. The RL algorithm applies to both discrete and continuous optimization parameters. The number of optimization variables is such that the  mixed combinatorial and continuous nature of the problem does not overwhelm the RL algorithm. In our case, the RL approach proved to be infeasible due to the large number of discrete choices. A system engineering approach to design of buildings is presented in \cite{GEYER200912}. The authors introduce structural and parametric diagrams using UML/SysML \cite{OMGSysML} language to describe the components of the design, the parametric equations that constrain the design variables, and the objective functions used in a multi-objective optimization approach. The Pareto front is searched using a multi-objective genetic algorithm with a multiple elitist strategy. The strategy consists of a ranking scheme that
gives non-dominated designs the highest probability to be selected for crossover. While the building design problem does present some similarities with our problem, it differs in that our parametric constraints are dynamic, i.e., we deal with DAEs. Thus the complexity is higher even when considering a modest number of components.

\section{Conclusions}

In this paper, we presented an automated design process utilizing a bottom-up approach. The process begins with an initial possibly large topology of universal components that is iteratively refined until a sparse solution is found. The initial design is based on universal components, each of which can exhibit a range of behavior through basic components. This combination of modes and topology ensures a broad exploration of the design space. We demonstrated two approaches for addressing the combinatorial explosion typical of design optimization problems. The first approach relaxes discrete variables to continuous variables by transforming discrete switches into continuous switches. These continuous components are physically realizable, resulting in no loss in performance requirements. Additionally, sparsity is induced through an $L_1$ regularization cost that encourages the parameters of the continuous switches to be zero. The second approach employs a random search algorithm inspired in part by genetic algorithms. It combines discrete choices that select a topology with continuous optimization to compute optimal parameters for the chosen topology. The topology is iteratively simplified, or mutated, until no further improvement is obtained. Both approaches are supported by automated model simplification and reconstruction that reduce the complexity of the design model, in turn decreasing the time complexity for the continuous optimization algorithms that require model simulations. These continuous optimization algorithms are gradient-free. We are currently investigating the application of a differential programming paradigm to the design problem described in this paper, which would allow us to utilize gradient-based algorithms. The major challenge we face is extending automatic differentiation support to DAEs that typically require stiff, implicit numerical solvers.

\bibliographystyle{plain}
\bibliography{references,refs}

\end{document}